\documentclass[3p,11pt]{elsarticle}
\usepackage{color}
\usepackage[colorlinks,linkcolor=blue,citecolor=blue]{hyperref}
\usepackage{amsmath,amsfonts,amssymb,amsthm}
\usepackage{mathrsfs} 
\usepackage{graphicx,epstopdf,float,subfigure}
\usepackage[margin=10pt,font=small,labelfont=bf,labelsep=period]{caption}

\usepackage{multirow,booktabs}

\numberwithin{equation}{section} 
\newtheorem{theorem}{Theorem}[section]
\newtheorem{lemma}{Lemma}[section]

\newtheorem{example}{Example}[section]
\biboptions{sort&compress}

\newtheoremstyle{mythm}{1.5ex plus 1ex minus .2ex}{1.5ex plus 1ex minus .2ex}
{\song}{\parindent}{\song\bfseries}{}{1em}{}
\theoremstyle{mythm}

\begin{document}

\begin{frontmatter}

\title{A $C^0$-continuous nonconforming virtual element method for linear strain gradient elasticity}
\author{Jianguo Huang}
\ead{jghuang@sjtu.edu.cn}
\address{School of Mathematical Sciences, and MOE-LSC, Shanghai Jiao Tong University, Shanghai 200240, China}
\author{Yue Yu\footnote{Corresponding author.}}
\ead{terenceyuyue@xtu.edu.cn}
\address{Hunan Key Laboratory for Computation and Simulation in Science and Engineering, Key Laboratory of Intelligent Computing and Information Processing of Ministry of Education, National Center for Applied Mathematics in Hunan, School of Mathematics and Computational Science, Xiangtan University, Xiangtan, Hunan 411105, China}

\begin{abstract}
 A robust $C^0$-continuous nonconforming virtual element method (VEM) is developed for a boundary value problem arising from strain gradient elasticity in two dimensions, with the family of polygonal meshes satisfying a very general geometric assumption given in Brezzi et al. (2009) and Chen and Huang (2018). The stability condition of the VEMs is derived by establishing Korn-type inequalities and inverse inequalities. Some crucial commutative relations for locking-free analysis as in elastic problems are derived. The sharp
and uniform error estimates with respect to both the microscopic parameter and the Lam\'e coefficient
are achieved in the lowest-order case, which is also verified by numerical results.
\end{abstract}


\end{frontmatter}

\section{Introduction}

Elasticity problems have wide applications in mechanical, architectural, aerospace and materials engineering, as well
as in various fields of physics and life sciences. The mathematical theory of elasticity, especially the linearized theory, is one of the foundations of several engineering sciences, in which the classical continuum theories are frequently used for modelling, analyzing and predicting the behavior of macro-scale solids and structures. However, the ability of classical continuum theories for describing multi-scale phenomena is very limited due to the specific properties of materials at small scales. During the past decades, notably size effects of microstructural materials have been observed in various experimental studies as a result of the extensive applications of micro- and nano-scale detection equipment \cite{Fleck-Muller-Ashby-Hutchinson-1994,Stolken-Evans-1998,Poole-Ashby-Fleck-1996,Lam-Yang-Chongetal-2003,Tang-Alici-2011,Bazant-1999}.

Since the classical continuum theories cannot capture the scale-dependent phenomenon, scholars have put forward several improved models by introducing high order or generalized continuum theories to analyze the size-dependent behavior of small-scale structures. Different from the classical continuum mechanics theories, internal length scale parameters characteristic of the underlying microstructure are involved in the constitutive equations \cite{Toupin-1962,Mindlin-Tiersten-1962,Koiter-1964,Mindlin-1964,Askes-Aifantis-2011,Aifantis-1984,Altan-Aifantis-1992,Ru-Aifantis-1993},
among which the Aifantis' strain gradient elasticity (SGE) model has become a new research focus in elastic and plastic problems in the analysis of microstructural materials. In this model, the classical constitutive equation is modified by incorporating the Laplacian of the strain or stress as (cf. \cite{Altan-Aifantis-1992,Askes-Aifantis-2011})
\begin{equation}\label{constitutiveSGE}
\tilde {\boldsymbol\sigma}_{ij} = \mathbb{C}_{ijkl}(\boldsymbol\varepsilon_{kl} - \iota^2\boldsymbol\varepsilon_{kl,mm})\quad \mbox{or}\quad \boldsymbol{\tilde \sigma} = \boldsymbol\sigma  - \iota^2\Delta \boldsymbol\sigma,
\end{equation}
where $\boldsymbol\varepsilon $ and $\boldsymbol\sigma$ are the strain and stress tensors in elasticity, $\tilde {\boldsymbol\sigma} $ is the modified stress tensor, and $\iota $ is the characteristic length or microscopic parameter of the underlying materials. Moreover, for linear elasticity problem, $\mathbb{C}_{ijkl} = \lambda \delta_{ij}\delta_{kl} + 2\mu \delta_{ik}\delta_{jl} $ with $\delta$ denoting the Kronechker delta. Throughout the paper, we use the summation convention whereby summation is implied when an index is repeated exactly two times.

Within the classical elasticity theory, the linear elastic model leads to a second order partial differential equation. For gradient elastic models, instead, a fourth-order governing equation will be induced.
For this reason, the finite element methods for plate bending problem are widely applied in the engineering to numerically solve the gradient elastic problems, but most of them lack theoretical analysis.
Conforming finite element methods can be referred to \cite{Manzari-Yonten-2013,Papanicolopulos2009,Torabi-Ansari-Darvizeh-2018,Babu-Patel-2019,Borst-Pamin-1996,Niiranen-Khakalo-Balobanov-Niemi-2016,Song-Zhao-He-Zheng-2014}. Mixed finite element can be found in \cite{Shu-King-Fleck-1999,Amanatidou-Aravas-2002,Borst-Pamin-1996,Phunpeng-Baiz-2015,Liao-Ming-Xu-2021}.
An alternative approach is to use the nonconforming finite element \cite{Soh-Chen-2004,Zhao-Chen-Lo-2011,Zybell-Muhlich-Kuna-Zhang-2012,Torabi-Ansari-Darvizeh-2019}.
In recent years, based on the finite elements proposed in \cite{Mardal-Tai-Winther-1991,Nilssen-Tai-Winther-2000} to deal with the fourth-order singular perturbation problems, Ming and his coauthors proposed several robust nonconforming finite element methods to solve the SGE problem, where a systematic theoretical analysis is also carried out \cite{Li-Ming-Shi-2017,Li-Ming-Shi-2018,Liao-Ming-Xu-2021}. In particular, Ref. \cite{Liao-Ming-Xu-2021} presented the mixed formulation of the SGE model and proposed Taylor-Hood finite elements to resolve the nearly incompressible strain gradient elasticity. Based on this mixed formulation, Ref. \cite{CHH2022} rigorously established the sharp and uniform error estimates with respect to both the microscopic parameter and the Lam\'e coefficient, and developed a double-parameter robust nonconforming mixed finite element method.

In this paper, we are concerned with the virtual element methods. To the best of our knowledge, no related work has been carried out.
As a generalization of the standard finite element method that allows for general polytopal meshes, the VEM is first proposed and analyzed in \cite{Beirao-Brezzi-Cangiani-2013}. The other pioneering works can be found in \cite{Ahmad-Alsaedi-Brezzi-2013,Beirao-Brezzi-Marini-2014}. VEMs have some advantages over standard finite element methods. For example, they are more convenient to handle PDEs on complex geometric domains or the ones associated with high-regularity admissible spaces. Until now, there have developed many conforming and nonconforming VEMs for second order elliptic equations \cite{Ahmad-Alsaedi-Brezzi-2013,Beirao-Brezzi-Cangiani-2013,DeDios-Lipnikov-Manzini-2016,Cangiani-Manzini-Sutton-2016} and fourth order elliptic equations \cite{Brezzi-Marini-2013,Chinosi-Marini-2016,Antonietti-Manzini-Verani-2018,Zhao-Zhang-Chen-2018}, respectively. From our standpoints, there are two features worth mention for VEMs. The first one is that the construction of nonconforming VEMs for elliptic problems is very natural and standard, which can be derived based on an integration by parts formula for the underlying differential operator. As a matter of fact, this idea is used to devise nonconforming VEMs for arbitrary order elliptic problems though the resulting formulation and theoretical analysis are rather involved \cite{Chen-HuangX-2020}. The other one is that if polytopal meshes reduce to simplicial ones, VEMs often give rise to some standard finite elements directly.

To avoid complicated presentation, we confine our discussion in two dimensions, with the family of polygonal meshes  $\{ \mathcal{T}_h \}_h$ satisfying the following condition (cf. \cite{Brezzi-Buffa-Lipnikov-2009,Chen-HuangJ-2018}):
\begin{enumerate}
\item[{\bf A1}.] For each $K\in {\mathcal T}_h$, there exists a ``virtual triangulation" ${\mathcal T}_K$ of $K$ such that ${\mathcal T}_K$ is uniformly shape regular and quasi-uniform. The corresponding mesh size of ${\mathcal T}_K$ is uniformly proportional to $h_K$. Each edge of $K$ is a side of a certain triangle in ${\mathcal T}_K$.
\end{enumerate}
As shown in \cite{Chen-HuangJ-2018}, this condition covers the usual geometric assumptions frequently used in the context of VEMs.

We propose and analyze the $C^0$-continuous nonconforming VEMs for the strain gradient elasticity.
Motivated by some ideas in \cite{Ahmad-Alsaedi-Brezzi-2013,Beirao-Brezzi-Marini-2013,Brezzi-Marini-2013,Zhao-Chen-Zhang-2016,Zhang-Zhao-Yang-2019,Zhang-Zhao-Chen-2020}, we first design the $C^0$-continuous nonconforming VEMs in Section \ref{sec:c0VEM} by establishing the integration by parts formulas of the involved operators and analyzing the computability.
As the preparation for the error analysis, some estimates for a $C^0$-continuous nonconforming virtual element for fourth-order problems are investigated in Section \ref{sec:C0estimates}, including the inverse inequalities, the norm equivalence and the interpolation error estimate. It is worth pointing out that the above work has completely supplemented the missing analysis in \cite{Zhao-Chen-Zhang-2016} for the $C^0$-continuous virtual element method for the plate bending problem.

In view of the rephrased stability conditions \eqref{mod1} and \eqref{mod2}, we further establish Korn-type inequalities in Section \ref{sect:Korn-Norm} by employing a tricky argument, where the geometric dependence of the hidden constants is clearly presented under the assumption {\bf A1}. The combination of these ingredients lead naturally to the stability estimates of the underlying bilinear forms as required.
We derive the crucial commutative relations in the virtual element spaces for locking-free analysis as in elastic problems and establish the Strang-type lemma with the hidden constant independent of the model parameters. We hence deduce the robustness with respect to the Lam\'{e} constant and the microscopic parameter in the lowest-order case by using the regularity estimate established in \cite{CHH2022}.

We end this section by introducing some notations and symbols frequently used in this paper. For a bounded Lipschitz domain $D$, the symbol $( \cdot , \cdot )_D$ denotes the $L^2$-inner product on $D$, $\|\cdot\|_{0,D}$ denotes the $L^2$-norm, and $|\cdot|_{s,D}$ is the $H^s(D)$-seminorm. For all integer $k\ge 0$, $\mathbb{P}_k(D)$ is the set of polynomials of degree $\le k$ on $D$.
Let $e \subset \partial K$ be the common edge for elements $K = K^-$ and $K^+$, and let $v$ be a scalar function defined on $e$. We introduce the jump of $v$ on $e$ by $[v] = v^- - v^+$, where $v^-$ and $v^+$ are the traces of $v$ on $e$ from the interior and exterior of $K$, respectively. Moreover, for any two quantities $a$ and $b$, ``$a\lesssim b$" indicates ``$a\le C b$" with the constant $C$ independent of the mesh size $h_K$, and ``$a\eqsim b$" abbreviates ``$a\lesssim b\lesssim a$".
To present the degrees of freedom (d.o.f.s), we introduce a scaled monomial $\mathbb{M}_r(D)$ on a $d$-dimensional domain $D$
\[
\mathbb  M_{r} (D):= \Big \{ \Big ( \frac{\boldsymbol x -  \boldsymbol x_D}{h_D}\Big )^{\boldsymbol  s}, \quad |\boldsymbol  s|\le r\Big \},
\]
where $h_D$ is the diameter of $D$, $\boldsymbol  x_D$ is the centroid of $D$, and $r$ is a non-negative integer. For the multi-index ${\boldsymbol{s}} \in {\mathbb{N}^d}$, we follow the usual notation $\boldsymbol{x}^{\boldsymbol{s}} = x_1^{s_1} \cdots x_d^{s_d}$, where $|\boldsymbol{s}| = s_1 +  \cdots  + s_d$. Conventionally, $\mathbb  M_r (D) =\{0\}$ for $r\le -1$.

We also refer to \cite{Huang-Yu-2021} for some basic estimates in VEM analysis.

\section{Strain gradient elastic model} \label{sect:sgemodel}

The strain gradient elastic model is described by the following boundary value problem: Find $\boldsymbol u$ the displacement vector that solves (cf. \cite{Li-Ming-Shi-2017,Altan-Aifantis-1992,Askes-Aifantis-2011})
 \begin{equation}\label{Pro:SGE}
 \begin{cases}
  - {\rm div}\boldsymbol{\tilde \sigma }(\boldsymbol{u}) = \boldsymbol{f} \quad & \mbox{in}~\Omega,  \\
   \boldsymbol{u}  = \partial_{\boldsymbol n} \boldsymbol{u}  = \boldsymbol 0 \quad & \mbox{on}~\partial \Omega,
\end{cases}
\end{equation}
where $\boldsymbol{\tilde \sigma}(\boldsymbol u)$ is given in \eqref{constitutiveSGE}.
The variational problem of \eqref{Pro:SGE} is to find $\boldsymbol u \in \boldsymbol H_0^2(\Omega )$ such that
\begin{equation}\label{ContinuousP}
a(\boldsymbol u,\boldsymbol v) = (\boldsymbol f,\boldsymbol v),\quad \boldsymbol v \in \boldsymbol{H}_0^2(\Omega ),
\end{equation}
where
$a(\boldsymbol u,\boldsymbol v)
= ( \mathbb{C}\boldsymbol\varepsilon (\boldsymbol u),\boldsymbol\varepsilon(\boldsymbol v) )
+ \iota^2( \mathbb{D}\nabla \boldsymbol\varepsilon(\boldsymbol u),\nabla \boldsymbol\varepsilon(\boldsymbol v) )$,
and $\mathbb{D}$ is a six-order tensor defined by $\mathbb{D}_{ijklmn} = \lambda \delta_{il}\delta_{jk}\delta_{mn} + 2\mu \delta_{il}\delta_{jm}\delta_{kn}$.
Let
$c^K(\boldsymbol u,\boldsymbol v) = ( \mathbb{C}\boldsymbol\varepsilon (\boldsymbol u),\boldsymbol\varepsilon (\boldsymbol v) )_K$ and $d^K(\boldsymbol u,\boldsymbol v) = ( \mathbb{D}\nabla \boldsymbol\varepsilon (\boldsymbol u),\nabla \boldsymbol\varepsilon (\boldsymbol v) )_K$.
One can check that for $\boldsymbol v \in \boldsymbol{H}_0^2(\Omega )$ they can further split as
\[c^K(\boldsymbol u,\boldsymbol v) = 2\mu c_1^K(\boldsymbol u, \boldsymbol v) + \lambda c_2^K(\boldsymbol u, \boldsymbol v),\qquad
d^K(\boldsymbol u,\boldsymbol v) = 2\mu d_1^K(\boldsymbol u, \boldsymbol v) + \lambda d_2^K(\boldsymbol u, \boldsymbol v),\]
where
\[c_1^K(\boldsymbol u,\boldsymbol v) = ( \boldsymbol\varepsilon(\boldsymbol u),\boldsymbol\varepsilon(\boldsymbol v))_K, \quad
c_2^K(\boldsymbol u,\boldsymbol v) = ({\rm div}\boldsymbol u,{\rm div}\boldsymbol v)_K,\]
\[d_1^K(\boldsymbol u,\boldsymbol v) = (\nabla\boldsymbol\varepsilon (\boldsymbol u),\nabla\boldsymbol\varepsilon (\boldsymbol v))_K, \quad
d_2^K(\boldsymbol u,\boldsymbol v) = ( \nabla{\rm div}\boldsymbol u,\nabla{\rm div}\boldsymbol v)_K.\]
A detailed calculation for above equations can be found in Lemma \ref{lem:Eh}.
In what follows, we call the first bilinear function as the linear elastic part and the second one as the strain gradient elastic part, respectively.

We define an energy norm $\|\cdot\|_{\iota,h}$ by $\| \boldsymbol v \|_{\iota,h} = (| \boldsymbol v |_{1,h}^2 + \iota^2| \boldsymbol v |_{2,h}^2)^{1/2}$.
According to the classic Korn's inequality, there holds the coercivity (see Theorem 1 in \cite{Li-Ming-Shi-2017})
\[C(\Omega )\| \boldsymbol v \|_{\iota,h}^2 \le a(\boldsymbol v,\boldsymbol v) \le 2(\lambda  + \mu )\| \boldsymbol v \|_{\iota,h}^2,  \quad \boldsymbol v \in \boldsymbol H_0^2(\Omega).\]

\section{A $C^0$-continuous virtual element method} \label{sec:c0VEM}

\subsection{Green's formulas and the induced projectors} \label{Green}

In what follows, $\boldsymbol V_h(K)$ is supposed to be the desired local virtual element space. The precise definition will be given in the next subsection.

We first consider the linear elastic part.
Let $\boldsymbol u\in \boldsymbol H^2(K)$ and $\boldsymbol p \in (\mathbb{P}_k(K))^2$. Referring to Eq. (2.31) in \cite{Beirao-Brezzi-Marini-2013}, we have
 \begin{equation}\label{c1K}
 c_1^K(\boldsymbol u,\boldsymbol p) = \int_K \boldsymbol\varepsilon (\boldsymbol u):\boldsymbol\varepsilon (\boldsymbol p) {\rm d}x =  - \int_K {\rm div}(\boldsymbol\varepsilon (\boldsymbol p)) \cdot \boldsymbol u {\rm d}x + \int_{\partial K} (\boldsymbol\varepsilon (\boldsymbol p)  \boldsymbol n) \cdot \boldsymbol u {\rm d}s.
 \end{equation}
The induced elliptic projector
$\Pi_K^1: \boldsymbol V_h(K) \to (\mathbb{P}_k(K))^2$, $\boldsymbol v \mapsto \Pi_K^1\boldsymbol v$
is defined by
\[c_1^K(\Pi_K^1\boldsymbol u,\boldsymbol p) = c_1^K(\boldsymbol u,\boldsymbol p),~~\boldsymbol p \in (\mathbb{P}_k(K))^2.\]
Since $c_1^K(\cdot, \cdot)$ is only semi-positive definite, we impose the constraint
\begin{equation}\label{constrainPiC}
\int_{\partial K} (\Pi_K^1\boldsymbol u) \cdot \boldsymbol p {\rm d}s = \int_{\partial K} \boldsymbol u \cdot \boldsymbol p {\rm d}s,\quad \boldsymbol p \in {\boldsymbol{\rm RM}}(K) = \{ \boldsymbol v \in \boldsymbol{H}^1(K ): \boldsymbol\varepsilon (\boldsymbol v) = O \}.
\end{equation}

For $c_2^K$, we note that for every $\boldsymbol v\in \boldsymbol H^1(K)$ and for every $p\in \mathbb{P}_k(K)$, there holds
\[\int_K ({\rm div} \boldsymbol v)p {\rm d}x =  - \int_K {\boldsymbol v \cdot \nabla p} {\rm d}x + \int_{\partial K} (\boldsymbol v \cdot \boldsymbol n)p {\rm d}s.\]
Accordingly, we can introduce an ``$L^2$ projector'' using external discretization defined by (see Eq. (3.6) in \cite{Beirao-Brezzi-Marini-2013})
\begin{align}
  \Pi_{k - 1}^0{\rm div}: \boldsymbol V_h(K) \to \mathbb{P}_{k-1}(K), \quad \boldsymbol v \mapsto \Pi_{k - 1}^0({\rm div}\boldsymbol v), \nonumber\\
  \int_K \Pi_{k - 1}^0({\rm div} \boldsymbol v)p {\rm d}x = \int_K ({\rm div}\boldsymbol v)p {\rm d}x,~~p \in \mathbb{P}_{k - 1}(K). \label{L2ProC}
\end{align}

We next focus on the strain gradient elastic part. The Green's formula in the following lemma is used to define the elliptic projector $\Pi_K^2$.
\begin{lemma}\label{lem:d1k}
For all $\boldsymbol u \in \boldsymbol H^2(K)$ and for all $\boldsymbol p \in (\mathbb{P}_k(K))^2$, there holds the following Green's formula
\begin{align*}
  d_1^K(\boldsymbol u,\boldsymbol p) &= \int_K \boldsymbol Q_{3i,i}(\boldsymbol p) \cdot \boldsymbol u {\rm d}x - \int_{\partial K} ( \boldsymbol Q_{3 \boldsymbol n}(\boldsymbol p) + \partial_{\boldsymbol t}\boldsymbol M_{\boldsymbol{tn}}(\boldsymbol p) ) \cdot \boldsymbol u {\rm d}s  \\
   &\quad + \int_{\partial K} \boldsymbol M_{\boldsymbol{nn}}(\boldsymbol p) \cdot \partial_{\boldsymbol n}\boldsymbol u {\rm d}s - \sum\limits_{i = 1}^{N_v} \left[ \boldsymbol M_{\boldsymbol{tn}}(\boldsymbol p) \right](z_i) \cdot \boldsymbol u(z_i),
\end{align*}
where $[ M_{\boldsymbol{tn}}^k(\boldsymbol p) ](z_i) = M_{\boldsymbol{tn}}^k(\boldsymbol p)  |_{z_i^ -}^{z_i^ +}$ is the jump at the vertex $z_i$ along the boundary of $K$,
\[ \boldsymbol Q_{3i,i}^\intercal = [Q_{3i,i}^1,Q_{3i,i}^2], \boldsymbol Q_{3 \boldsymbol n}^\intercal = [Q_{3 \boldsymbol n}^1,Q_{3 \boldsymbol n}^2],
\boldsymbol M_{\boldsymbol{tn}}^\intercal = [M_{\boldsymbol{tn}}^1,M_{\boldsymbol{tn}}^2], \boldsymbol M_{\boldsymbol{nn}}^\intercal = [M_{\boldsymbol{nn}}^1,M_{\boldsymbol{nn}}^2],\]
and
\[M_{ij}^k(\boldsymbol p) = \partial_i\varepsilon_{kj}(\boldsymbol p),~~M_{\boldsymbol{nn}}^k(\boldsymbol p) = M_{ij}^k(\boldsymbol p)n_in_j,~~M_{\boldsymbol{tn}}^k(\boldsymbol p) = M_{ij}^k(\boldsymbol p)t_in_j,\quad k = 1,2,\]
\[Q_{3i}^k(\boldsymbol p) = M_{ij,j}^k(\boldsymbol p) = \partial_{ij}\varepsilon_{kj}(\boldsymbol p),\quad Q_{3 \boldsymbol n}^k(\boldsymbol p) = Q_{3i}^k(\boldsymbol p)n_i,\quad k = 1,2.\]
\end{lemma}
\begin{proof}
The formula follows from the integration by parts.
\end{proof}

The corresponding elliptic projection operator
\[\Pi_K^2: \boldsymbol V_h(K) \to (\mathbb{P}_k(K))^2, \quad \boldsymbol u \mapsto \Pi_K^2 \boldsymbol u\]
is defined by
$d_1^K(\Pi_K^2\boldsymbol u,\boldsymbol p) = d_1^K(\boldsymbol u,\boldsymbol p)$, $\boldsymbol p \in (\mathbb{P}_k(K))^2$.
Let $d_1^K(\boldsymbol p,\boldsymbol p) = 0$. One easily finds that $\boldsymbol p \in (\mathbb{P}_1(K))^2$. We impose the constraints
\begin{equation}\label{const2}
\int_{\partial K} {\Pi_K^2 \boldsymbol u} {\rm d}s = \int_{\partial K} \boldsymbol u {\rm d}s,\quad \int_{\partial K} {\nabla \Pi_K^2\boldsymbol u} {\rm d}s = \int_{\partial K} {\nabla \boldsymbol u} {\rm d}s.
\end{equation}
For the computation of the second condition, please refer to \cite{Zhao-Chen-Zhang-2016} or \eqref{compugrad}.

For $d_2^k$, following the idea in the linear elasticity problem, we are able to define the $L^2$ projection by considering $( {\nabla {\rm div}\boldsymbol u,\boldsymbol p} )_K$.
\begin{lemma}\label{lem:d2k}
For every $\boldsymbol u \in \boldsymbol H^2(K)$ and for every $\boldsymbol p \in (\mathbb{P}_{k - 2}(K))^2$, there holds
\begin{align*}
  ( \nabla {\rm div}\boldsymbol{u},\boldsymbol{p} )_K &= \int_K \boldsymbol{Q}_{1i,i}(\boldsymbol{p}) \cdot \boldsymbol{u} {\rm d}x - \int_{\partial K} ( \boldsymbol{Q}_{1\boldsymbol{n}}(\boldsymbol{p}) + \partial_{\boldsymbol{t}}\boldsymbol{M}_{0\boldsymbol{t}}(\boldsymbol{p}) ) \cdot \boldsymbol{u} {\rm d}x  \\
  &\quad + \int_{\partial K} \boldsymbol{M}_{0\boldsymbol{n}}(\boldsymbol{p}) \cdot \partial_{\boldsymbol{n}}\boldsymbol{u} {\rm d}s- \sum\limits_{i = 1}^{N_v} [ \boldsymbol{M}_{0\boldsymbol{t}}(\boldsymbol{p}) ](z_i) \cdot \boldsymbol{u}(z_i).
\end{align*}
where
\[\boldsymbol{Q}_{1i,i}^\intercal = [Q_{1i,i}^1,Q_{1i,i}^2], \boldsymbol{Q}_{1\boldsymbol n}^\intercal = [Q_{1\boldsymbol n}^1,Q_{1\boldsymbol n}^2],
\boldsymbol{M}_{0\boldsymbol t}^\intercal = {[M_{0\boldsymbol t}^1,M_{0\boldsymbol t}^2]},  \boldsymbol{M}_{0\boldsymbol n}^\intercal = [M_{0\boldsymbol{n}}^1,M_{0\boldsymbol{n}}^2],\]
\[M_{0\boldsymbol n}^j(\boldsymbol p) = {p_i}n_in_j,~~M_{0\boldsymbol t}^j(\boldsymbol p) = {p_i}t_in_j,~~Q_{1i}^j(\boldsymbol p) = {\partial _j}{p_i},~~Q_{1\boldsymbol n}^j(\boldsymbol p) = Q_{1i}^j(\boldsymbol p)n_i.\]
\end{lemma}

The induced ``$L^2$ projector''
\[\Pi _{k - 2}^0\nabla {\rm div}: \boldsymbol V_h(K) \to (\mathbb{P}_{k - 2}(K))^2,~~\boldsymbol u \mapsto \Pi _{k - 2}^0\nabla {\rm div}\boldsymbol u\]
is defined by $( \Pi _{k - 2}^0\nabla {\rm div} \boldsymbol u,\boldsymbol p )_K = ( \nabla {\rm div}\boldsymbol u,\boldsymbol p )_K,~ \boldsymbol p \in (\mathbb{P}_{k - 2}(K))^2$.

\subsection{The nonconforming virtual element space}

In \cite{Zhao-Chen-Zhang-2016}, a $C^0$-continuous nonconforming virtual element method has been proposed for plate bending problem, where the generic element $K$ is required to be convex. For this reason, we assume that every element is a convex polygon and always set $k\ge 2$. The proposed local virtual element space is defined by
\begin{align}
 V_k(K) =  \{ v \in H^2(K): \Delta^2v \in \mathbb{P}_{k - 4}(K),  v|_e \in \mathbb{P}_k(e),
\Delta v|_e \in \mathbb{P}_{k - 2}(e), e \subset \partial K  \},  \label{VkK}
\end{align}
with the usual convention that $\mathbb{P}_{ - 1} = \mathbb{P}_{ - 2} = \{ 0\} $.
Moreover, the space can be equipped with the following degrees of freedom:
\begin{itemize}
  \item $\boldsymbol \chi_1(v)$: the values at the vertices of $K$,
  \[\chi_i(v) = v(z_i),\quad i = 1, \cdots ,{N_v}.\]
  \item $\boldsymbol \chi_2(v)$: the moments of $v$ on edges up to degree $k-2$,
  \[\chi_e(v) = |e|^{-1}(m_e,v)_e,\quad m_e \in \mathbb{M}_{k - 2}(e).\]
  \item $\boldsymbol \chi_3(v)$: the moments of $\partial_{\boldsymbol n}v$ on edges up to degree $k-2$,
  \[\chi_{n_e}(v) = (m_e,\partial_{\boldsymbol n}v)_e,\quad m_e \in \mathbb{M}_{k - 2}(e).\]
  \item $\boldsymbol \chi_4(v)$: the moments on element $K$ up to degree $k-4$,
  \[\chi_K(v) = |K|^{-1}(m_K,v),\quad m_K \in \mathbb{M}_{k - 4}(K).\]
\end{itemize}

Observing the Green's formulas in Lemmas \ref{lem:d1k} and \ref{lem:d2k}, one can check that the projections $\Pi_K^2 \boldsymbol v$ and $\Pi_{k - 2}^0\nabla \mathrm{div}\boldsymbol{v}$ are uniquely determined by the d.o.f.s of each component of $\boldsymbol v$, where $\boldsymbol v$ is a function in the tensor product space $(V_k(K))^2$.
We point out that, however, additional moments on element $K$ with degrees $k-3$ and $k-2$ are required for the computation of the projectors $\Pi_K^1$ and $\Pi_{k - 1}^0{\rm div}$ associated with the linear elastic part. For this reason, we enlarge $V_k(K)$ to be a lifting space as
\begin{equation*}\label{liftedSpace}
\widetilde V_k(K) = \left\{ v \in H^2(K):\Delta^2v \in \mathbb{P}_{k - 2}(K),~v|_e \in \mathbb{P}_k(e),~\Delta v|_e \in \mathbb{P}_{k - 2}(e),~e \subset \partial K \right\},
\end{equation*}
in which the functions are uniquely determined by the previous d.o.f.s together with the additional moments
$\chi_K(v) = |K|^{-1} (m, v)_K$ with $ m \in \mathbb{M}_{k - 2}(K)\backslash \mathbb{M}_{k - 4}(K)$.
We now define an elliptic projector on the lifting tensor product space as
\begin{align*}
  \widetilde \Pi_K^2: (\widetilde{V}_h(K))^2 \to (\mathbb{P}_k(K))^2,\quad \boldsymbol{v} \mapsto \widetilde \Pi_K^2\boldsymbol{v}, \\
  d_1^K(\widetilde \Pi_K^2\boldsymbol{v},\boldsymbol p) = d_1^K(\boldsymbol{v},\boldsymbol p),\quad \boldsymbol p \in (\mathbb{P}_k(K))^2,
\end{align*}
with the same constraints as $\Pi_K^2$.

In what follows, denote by $\boldsymbol{V}_h = (\widetilde{V}_h)^2$ the global lifting space consisting of nonconforming virtual elements of the same type as mentioned above, which is defined elementwise and required that the degrees of freedom are continuous through the interior vertices and edges while zero for the boundary d.o.f.s.

\subsection{Construction of the discrete problem}

With the previous projections in mind, we are in a position to construct the approximate bilinear form.
The discrete variational problem is: Find $\boldsymbol u_h \in \boldsymbol V_h$ such that
\begin{equation}\label{DiscreteP}
a_h(\boldsymbol u_h,\boldsymbol v_h) = \langle \boldsymbol f_h,\boldsymbol v_h\rangle ,\quad \boldsymbol v_h \in \boldsymbol V_h,
\end{equation}
where
\begin{equation}\label{ahK}
a_h^K(\boldsymbol u,\boldsymbol v) = c_h^K(\boldsymbol u,\boldsymbol v) + \iota^2d_h^K(\boldsymbol u,\boldsymbol v).
\end{equation}
The first bilinear form on the right-hand side is split as
$c_h^K(\boldsymbol u,\boldsymbol v) = 2\mu c_{1,h}^K(\boldsymbol u,\boldsymbol v) + \lambda c_{2,h}^K(\boldsymbol u,\boldsymbol v)$,
where
\begin{align*}
  c_{1,h}^K(\boldsymbol u,\boldsymbol v) & = c_1^K(\Pi_K^1\boldsymbol u,\Pi_K^1\boldsymbol v) + S^K(\boldsymbol u - \Pi_K^1\boldsymbol u,\boldsymbol v - \Pi_K^1\boldsymbol v), \\
  c_{2,h}^K(\boldsymbol u,\boldsymbol v) & = (\Pi_{k - 1}^0{\rm div}\boldsymbol u,\Pi_{k - 1}^0{\rm div}\boldsymbol v)_K.
\end{align*}
The second one associated with the strain gradient elastic part is decomposed as
$d_h^K(\boldsymbol u,\boldsymbol v) = 2\mu d_{1,h}^K(\boldsymbol u,\boldsymbol v) + \lambda d_{2,h}^K(\boldsymbol u,\boldsymbol v)$,
where
\begin{align*}
  d_{1,h}^K(\boldsymbol u,\boldsymbol v) & = d_1^K(\Pi_K^2\boldsymbol u,\Pi_K^2\boldsymbol v) + h_K^{ - 2}S^K(\boldsymbol u - \Pi_K^2\boldsymbol u,\boldsymbol v - \Pi_K^2\boldsymbol v), \\
  d_{2,h}^K(\boldsymbol u,\boldsymbol v) & = (\nabla \Pi_{k - 1}^0{\rm div}\boldsymbol u,\nabla \Pi_{k - 1}^0{\rm div}\boldsymbol v)_K.
\end{align*}

A computable stabilization term $S^K(\cdot, \cdot)$ is constructed such that there holds
\begin{itemize}
  \item $k$-consistency: If $\boldsymbol v \in \boldsymbol V_h(K)$ and $\boldsymbol p \in (\mathbb{P}_k(K))^2$, then
\begin{equation}\label{consistency}
a_h^K(\boldsymbol v,\boldsymbol p)=a_K(\boldsymbol v,\boldsymbol p).
\end{equation}
  \item Stability: There exist constants ${\alpha_*}$ and ${\alpha^*}$, independent of $h$, $\iota$ and $\lambda$, such that
\begin{align}
  \alpha_*c_1^K(\boldsymbol v,\boldsymbol v) \le c_{1,h}^K(\boldsymbol v,\boldsymbol v) \le \alpha^*c_1^K(\boldsymbol v,\boldsymbol v),\quad \boldsymbol v \in \boldsymbol V_h(K), \label{stabilityC} \\
  \alpha_*d_1^K(\boldsymbol v,\boldsymbol v) \le d_{1,h}^K(\boldsymbol v,\boldsymbol v) \le \alpha^*d_1^K(\boldsymbol v,\boldsymbol v),\quad \boldsymbol v \in \boldsymbol V_h(K), \label{stabilityD}
\end{align}
\end{itemize}
Under the stability conditions, one easily finds that
\[\| \boldsymbol{\varepsilon}(\boldsymbol{v}) \|_{0,K}^2 + \iota^2\| \nabla \boldsymbol{\varepsilon}(\boldsymbol{v}) \|_{0,K}^2 \lesssim a_h^K(\boldsymbol{v},\boldsymbol{v}) ,\quad \boldsymbol{v} \in \boldsymbol{V}_h(K).\]
In view of the $C^0$ continuity of the functions in $\boldsymbol V_h$, the discrete Korn's inequality (1.15) in \cite{Brenner-2004} implies
\begin{equation}\label{KornOmega}
|\boldsymbol v_h|_{1,h} \lesssim \| \boldsymbol\varepsilon (\boldsymbol v_h)\|_{0,h}\quad \boldsymbol v_h\in \boldsymbol V_h.
\end{equation}
Observing the following equality
\begin{equation}\label{combieps}
  \frac{\partial^2v_i}{\partial x_j\partial x_k} = \frac{\partial}{\partial x_j}{\varepsilon_{ik}}(\boldsymbol v) + \frac{\partial}{\partial x_k}{\varepsilon_{ij}}(\boldsymbol v) - \frac{\partial}{\partial x_i}{\varepsilon_{jk}}(\boldsymbol v),\quad 1\le i,j,k\le 2,
\end{equation}
we immediately obtain
\begin{equation}\label{Korngrad}
\| \nabla \boldsymbol\varepsilon (\boldsymbol v) \|_{0,K} \eqsim | \boldsymbol v |_{2,K},\quad \boldsymbol v\in \boldsymbol H^2(K),
\end{equation}
which along with \eqref{KornOmega} yields the coercivity
\begin{equation}\label{coercivity}
\| \boldsymbol v_h \|_{\iota,h} \lesssim  a_h(\boldsymbol{v}_h,\boldsymbol{v}_h),\quad \boldsymbol{v}_h \in \boldsymbol{V}_h.
\end{equation}

It is evident that the stability conditions can be reformulated as
\begin{equation}\label{refstab1}
S^K(\boldsymbol{\tilde v},\boldsymbol{\tilde v}) \eqsim c_1^K(\boldsymbol{\tilde v},\boldsymbol{\tilde v}) = \| \boldsymbol\varepsilon (\boldsymbol{\tilde v}) \|_{0,K}^2,\quad \boldsymbol{\tilde v} = \boldsymbol v - \Pi_K^1\boldsymbol v
\end{equation}
and
\begin{equation}\label{refstab2}
h_K^{ - 2}S^K(\boldsymbol{\tilde w},\boldsymbol{\tilde w}) \eqsim d_1^K(\boldsymbol{\tilde w},\boldsymbol{\tilde w})=\| \nabla \boldsymbol\varepsilon (\boldsymbol{\tilde w}) \|_{0,K}^2,\quad \boldsymbol{\tilde w} = {\boldsymbol w} - \Pi_K^2{\boldsymbol w},
\end{equation}
which will be proved by establishing the following relations
\begin{align}
S^K(\boldsymbol{\tilde v},\boldsymbol{\tilde v}) & \eqsim | \boldsymbol{\tilde v} |_{1,K}^2 \eqsim \| \boldsymbol\varepsilon (\boldsymbol{\tilde v}) \|_{0,K}^2 ,\quad \boldsymbol{\tilde v} = \boldsymbol v - \Pi_K^1\boldsymbol v \label{mod1}\\
h_K^{-2}S^K(\boldsymbol{\tilde w},\boldsymbol{\tilde w}) & \eqsim | \boldsymbol{\tilde w} |_{2,K}^2 \eqsim \| \nabla \boldsymbol\varepsilon (\boldsymbol{\tilde w}) \|_{0,K}^2 ,\quad \boldsymbol{\tilde w} = \boldsymbol w - \Pi_K^2\boldsymbol w. \label{mod2}
\end{align}
The first ``$ \eqsim $'' and the second ``$ \eqsim $'' in \eqref{mod1} and \eqref{mod2} are related to the so-called norm equivalence and Korn's inequalities, respectively.

Relabel the d.o.f.s of $V_k(K)$ by a single index $i = 1, 2,  \cdots, N_K:= \dim V_k(K)$ to form a set $\{\chi_i\}_{i=1}^{N_K}$.
The stabilization term is then realized as
$S^K(\boldsymbol v,\boldsymbol w) = \sum\limits_{i = 1}^{N_K} \chi_i(\boldsymbol v)\cdot\chi_i(\boldsymbol w)$,
where $\chi_i(\boldsymbol v) = [\chi_i(v_1), \chi_i(v_2)]^\intercal$, $\boldsymbol v = [v_1,v_2]^\intercal$ and $\boldsymbol w = [w_1,w_2]^\intercal$.

For $k\ge 2$, the right-hand side can be simply approximated by (cf. \cite{Zhang-Zhao-Chen-2020})
\[\boldsymbol f_h |_K = \Pi_{k - 2}^0\boldsymbol f := [\Pi_{k - 2}^0f_1,\Pi_{k - 2}^0f_2]^\intercal,\quad K \in \mathcal{T}_h.\]


\section{Some estimates of a $C^0$-continuous and $H^2$-nonconforming virtual element} \label{sec:C0estimates}

Let $k,l$ be two nonnegative integers with $k\ge 2$. We introduce a virtual element space on $K$ as
\begin{equation}\label{Vkl}
V_{k,l}(K) = \left\{ v \in H^2(K): \Delta^2v \in \mathbb{P}_l(K), v|_e \in \mathbb{P}_k(e), \Delta v|_e \in \mathbb{P}_{k - 2}(e), e \subset \partial K \right\}.
\end{equation}
The d.o.f.s are given by
\begin{itemize}
  \item $\boldsymbol\chi_1(v)$: the values at the vertices of $K$;
  \item $\boldsymbol\chi_2(v)$: the moments of $v$ on edges up to degree $k-2$;
  \item $\boldsymbol\chi_3(v)$: the moments of $\partial_{\boldsymbol n}v$ on edges up to degree $k-2$;
  \item $\boldsymbol\chi_4(v)$: the moments on the element $K$ up to degree $l$.
\end{itemize}
Obviously, when $l=k-4$, $V_{k,l}(K)$ is the space $V_k(K)$ defined in \eqref{VkK}; when $l=k-2$, $V_{k,l}(K)$ is the lifting space $\widetilde{V}_k(K)$ given by \eqref{liftedSpace}. We in this section present some basic estimates for the above $C^0$-continuous nonconforming virtual element space $V_{k,l}(K)$.

We first consider the inverse inequalities, which play fundamental role in establishing the stability conditions. We will utilize the minimum energy principle based method and apply a modified idea in \cite{Brenner-Guan-Sung-2017} to construct the desired function. To this end, we define the following mesh-dependent norm
\begin{equation}\label{meshnorm}
|\!|\!|v|\!|\!|_K = \| \Pi_l^0v \|_{0,K} + h_K^{ 1/2}\| v \|_{0,\partial K} + h_K^{ 3/2}\Big( \sum\limits_{e \subset \partial K} \| \Pi_{k - 2,e}^0\partial_{\boldsymbol n}v \|_{0,e}^2 \Big)^{1/2},
\end{equation}
which plays the role of $L^2$ norm in VEM analysis and can be chosen as the candidate of stabilization terms (modified as a bilinear form). In fact, it is equivalent to the stabilization term induced by the d.o.f. vector.
\begin{lemma}\label{lem:normeqchi}
For all $v\in V_{k,l}(K)$ there hold the following norm equivalence formulas:
\begin{align}
  &h_K^{-1}( \| \boldsymbol\chi_1(v) \|_{l^2}^2 + \| \boldsymbol\chi_2(v) \|_{l^2}^2 )^{1/2} \eqsim h_K^{ - 3/2}\Big( \sum\limits_{e \subset \partial K}\|v\|_{0,e}^2  \Big)^{1/2}, \label{norm12} \\
  &h_K^{-1}\| \boldsymbol\chi_3(v) \|_{l^2} \eqsim h_K^{ - 1/2}\Big( \sum\limits_{e \subset \partial K} \| \Pi_{k - 2,e}^0\partial_{\boldsymbol n}v \|_{0,e}^2  \Big)^{1/2}, \label{norm3} \\
  &h_K^{-1}\|\boldsymbol\chi_4(v)\|_{l^2} \eqsim h_K^{-2}\| \Pi_l^0v \|_{0,K}. \label{norm4}
\end{align}
\end{lemma}
\begin{proof}
The first equation follows from the norm equivalence in finite element methods since $v|_e$ is a polynomial for $e\subset \partial K$.

For the third equation, let $m_\alpha$ be the scaled monomial on $K$ with degree no more than $l$. Let $\Pi_l^0v = \sum\nolimits_\alpha  c_\alpha m_\alpha$ and denote $\boldsymbol{c} = (c_\alpha)$ the coefficient vector. It is evident that $\|\Pi_l^0v\|_{0,K}^2 = \boldsymbol{c}^\intercal \boldsymbol{M} \boldsymbol{c} $, where $\boldsymbol{M} = ( (m_\alpha ,m_\beta )_K)$. According to the norm equivalence in  \cite{Chen-HuangJ-2018} (see Lemma 4.1 there), we have
$h_K^2\|\boldsymbol{c}\|_{l^2}^2 \eqsim \|\Pi_l^0v\|_{0,K}^2 = \boldsymbol{c}^\intercal \boldsymbol{M} \boldsymbol{c}$.
From the Rayleigh representation theorem of eigenvalues of symmetric positive definite matrices, one obtains that any eigenvalue of $\boldsymbol{M}$ satisfies $\lambda(\boldsymbol{M})\eqsim h_K^2$, which gives
$\|\boldsymbol{Mc}\|_{l^2} \eqsim h_K^2 \|\boldsymbol{c}\|_{l^2}$.
Taking the internal moments on the expansion of $\Pi_l^0v$ yields
$\boldsymbol\chi_4(\Pi_l^0v) = |K|^{-1} \boldsymbol{Mc}$.
According to the definition of $L^2$ projections, one has
$|K|^{-1}(\Pi_l^0v, m_K)_K = |K|^{-1}(v,m_K)_K,   m_K \in \mathbb{M}_l(K)$,
which means $\boldsymbol\chi_4(v) = \boldsymbol\chi_4(\Pi_l^0v)$. Combining above equations implies that
\[\|\boldsymbol\chi_4(v)\|_{l^2} = \|\boldsymbol\chi_4(\Pi_l^0v)\|_{l^2} = |K|^{-1} \|\boldsymbol{Mc}\|_{l^2} \eqsim \|\boldsymbol{c}\|_{l^2}\eqsim h_K^{-1} \|\Pi_l^0v\|_{0,K}.\]
This is exactly the third equation \eqref{norm4}.

For the second one, we only need to take $w=\partial_{\boldsymbol n}v$ and get the desired result by utilizing the similar argument above on each edge $e$.
\end{proof}

\begin{lemma}\label{rigA2}
 Let $K$ be an element and let $w$ be a function defined on $K$ which satisfies
$\Delta w \in \mathbb{P}_l(K)$ and $ w|_e \in \mathbb{P}_k(e)$ for $e\subset \partial K$,
where $l,k$ are integers and $\mathbb{P}_k = \{0\}$, $k<0$. Then there holds
 $\|w\|_{0,\partial K} \lesssim h_K^{-1/2} \|w\|_{0,K}$.
\end{lemma}


\begin{proof}
According to Lemma 3.4 in \cite{Huang-Yu-2021}, there exists a polynomials $p$ such that $\Delta w = \Delta p$ and
$|p|_{1,K} \lesssim h_K^{-1}\|w\|_{0,K}$.
Applying the trace inequality to get
\[\|p\|_{0,\partial K} \lesssim h_K^{1/2}|p|_{1,K} + h_K^{-1/2}\|p\|_{0,K} \lesssim h_K^{-1/2}\|w\|_{0,K}.\]
By the triangle inequality it suffices to prove the desired estimate for $w-p$.  In other words, without loss of generality, we can additionally assume that $\Delta w = 0$.

For a boundary edge $e\subset \partial K$, denote by $\tau$ the triangle in the virtual triangulation with $e$ as an edge. Let $\lambda_1$ and $\lambda_2$ be the barycentric coordinate functions of $\tau$ associated with the endpoints of $e$ and denote $b_e = \lambda_1\lambda_2$. Since $w$ is a harmonic function, $w\in H^1(\tau)$ when restricted on $\tau$. Hence
\[b_ew \in H^1(\tau), \quad b_ew|_{\partial \tau \backslash e} = 0, \quad b_ew|_e~\mbox{is a polynomial}.\]
Let $w_1$ be a function satisfying
\[\begin{cases}
\Delta w_1 = 0 \quad \mbox{in $\tau$}, \\
w_1 = b_ew \quad \mbox{on $\partial \tau$}.
\end{cases}
\]
We obtain from the trace inequality and the inverse inequalities of conforming virtual element functions that
\begin{equation}\label{rig1}
\|w_1\|_{0,e} \lesssim h_K^{1/2}|w_1|_{1,\tau} + h_K^{-1/2}\|w_1\|_{0,\tau} \lesssim h_K^{-1/2}\|w_1\|_{0,\tau}.
\end{equation}

Let $w_2 = b_ew-w_1$. Then
$w_2|_{\partial \tau} = 0,   \Delta w_2 = \Delta b_e w + 2\nabla b_e \cdot \nabla w\in H^{-1}(\tau)$.
Integrating by parts and using the Poincar\'{e}-Friedrichs inequality, we derive
\[|w_2|_{1,\tau}^2 = -(\Delta w_2, w_2)_\tau \le \|\Delta w_2\|_{-1,\tau} \|w_2\|_{1,\tau}\lesssim \|\Delta w_2\|_{-1,\tau} |w_2|_{1,\tau},\]
which gives
$|w_2|_{1,\tau} \lesssim \|\Delta w_2\|_{-1,\tau}$.
Applying the Poincar\'{e}-Friedrichs inequality once again yields
\begin{align*}
\|w_2\|_{0,\tau}
& \lesssim h_K|w_2|_{1,\tau} \lesssim h_K\|\Delta w_2\|_{-1,\tau}
 \le h_K(\|\Delta b_e w\|_{-1,\tau} + \| \nabla b_e \cdot \nabla w\|_{-1,\tau})\\
& \lesssim h_K( h_K^{-2}\|w\|_{-1,\tau} + h_K^{-1}\|w\|_{0,\tau} ).
\end{align*}
According to the definition of negative norms and the Poincar\'{e}-Friedrichs inequality, one has
\[\|w\|_{-1,\tau} = \sup_{q\in H_0^1(\tau)} \frac{(w,q)_\tau}{\|q\|_{1,\tau}}\lesssim h_K\|w\|_{0,\tau}.\]
Plugging it in the previous inequality, we then have
$\|w_2\|_{0,\tau}\lesssim \|w\|_{0,\tau}$
and hence get
\begin{equation}\label{rig2}
\|w_1\|_{0,\tau} \le \|b_ew\|_{0,\tau} + \|w_2\|_{0,\tau} \lesssim \|w\|_{0,\tau}.
\end{equation}

From the properties of bubble functions one has
$\|w\|_{0,e}\lesssim \|b_e^{1/2}w\|_{0,e}\le \|b_ew\|_{0,e} = \|w_1\|_{0,e}$,
which together with \eqref{rig1} and \eqref{rig2} implies the desired result.
\end{proof}

In view of the relations in Lemma \ref{lem:normeqchi}, we now establish the inverse inequality in $V_{k,l}(K)$ with respect to the mesh-dependent norm \eqref{meshnorm}.
\begin{lemma}\label{lem:invH2mesh}
  There holds
 $| v |_{2,K} \lesssim h_K^{ - 2}|\!|\!|v|\!|\!|_K$ for all $v \in V_{k,l}(K)$.
\end{lemma}
\begin{proof}
Lemma 2.19 in \cite{Brenner-Guan-Sung-2017} uses the minimum energy principle to prove the corresponding result of the virtual element space of the second-order problem. The core of the argument is to construct auxiliary functions using Riesz representation theorem. However, the extension to fourth-order problems may be not straightforward. We now present a new construction method. For clarity, we divide the proof into several steps.

Step 1: Observing the Green's formula
\begin{align*}
  \int_K \nabla^2v: \nabla^2 (v - w) {\rm d}x
  &= \int_K \Delta^2v(v - w) {\rm d}x - \int_{\partial K} \partial_{\boldsymbol n}(\Delta v)(v - w) {\rm d}s \\
  &\quad + \int_{\partial K} (\Delta v - \partial_{\boldsymbol t}^2v)\partial_{\boldsymbol n}(v - w) {\rm d}s + \int_{\partial K} \partial_{\boldsymbol{\boldsymbol{tn}}}^2v\partial_{\boldsymbol t}(v - w) {\rm d}s,
\end{align*}
we obtain
$\int_K \nabla^2v: \nabla^2(v - w) {\rm d}x = 0$ for $w\in S(K)$,
where $S(K)$ is a subspace of $H^2(K)$ and every $w\in S(K)$ satisfies
$w|_e = v|_e, \Pi_{k-2,e}^0(\partial_{\boldsymbol n}w|_e) = \Pi_{k-2,e}^0(\partial_{\boldsymbol n}v|_e)$
for every $e \subset \partial K$ and
$\Pi_{k - 4}^0w = \Pi_{k - 4}^0v$.
This implies the minimum energy principle
$| v |_{2,K} \le | w |_{2,K}$ for all $ w \in S(K)$.
We will construct a function $w:=v_K \in S(K)$ to yield the desired estimate.

Step 2: For the virtual triangulation $\mathcal{T}_K$ of $K$, on each triangle we employ a $\mathbb{P}_k$ macroelement defined in \cite{Douglas-Dupont-Percell-1979}.  We choose $S_h(K)$ to be a subspace of the macroelement space associated with $\mathcal{T}_K$ with vanishing d.o.f.s except that $\zeta_0 = v$ on $\partial K$. It is evident that $\zeta_0\in H^2(K)$. According to the mesh assumption {\bf A1}, we obtain from the inverse inequalities for polynomials and scaling arguments that
\begin{equation}\label{scalingzeta0}
h_K^2| \zeta_0|_{2,K} \lesssim h_K|\zeta_0|_{1,K} \lesssim \|\zeta_0\|_{0,K} \eqsim h_K^{1/2}\|\zeta_0\|_{0,\partial K} = h_K^{1/2}\| v \|_{0,\partial K},
\end{equation}
which gives
\begin{equation}\label{zeta0L2}
\| \zeta_0\|_{0,K} \lesssim |\!|\!|v|\!|\!|_K, \quad | \zeta_0|_{2,K} \lesssim h_K^{ - 2}|\!|\!|v|\!|\!|_K.
\end{equation}

Step 3: Let $\zeta  = \zeta_0 + \zeta_1 $, where $\zeta_1\in H^2(K)$ is a $C^0$-continuous nonconforming virtual element function satisfying
\begin{equation}\label{zeta1construct}
\Delta^2\zeta_1 = 0, \quad  \zeta_1|_e = 0, \quad \Delta \zeta_1|_e\in \mathbb{P}_{k-2}(e), \quad e\subset \partial K.
\end{equation}
This means the following d.o.f.s of $\zeta_1$ are zero: values at vertices, moments of $\zeta_1$ on edges up to $k-2$ and moments up to $l$ on $K$. We additionally require that
\[\Pi_{k - 2,e}^0\partial_{\boldsymbol n}\zeta  = \Pi_{k - 2,e}^0\partial_{\boldsymbol n}v \quad \mbox{or} \quad
\Pi_{k - 2,e}^0\partial_{\boldsymbol n}\zeta_1  = \Pi_{k - 2,e}^0\partial_{\boldsymbol n}(v-\zeta_0), \quad e\subset \partial K,\]
which determines the remaining d.o.f.s of $\zeta_1$: moments of $\partial_{\boldsymbol n}\zeta_1$ up to $e\subset \partial K$. It is obvious that $\zeta_1\in H^2(K)\cap H_0^1(K)$ and $\zeta|_e = v|_e$ for $e\subset \partial K$.

Using the integration by parts and the Cauchy-Schwarz inequality we obtain
\begin{align}
  |\zeta_1|_{2,K}^2
  & = \int_{\partial K} \Delta \zeta_1 \partial_{\boldsymbol n}\zeta_1 {\rm d}s
    = \sum\limits_{e\subset \partial K}\int_e \Delta \zeta_1 \Pi_{k-2,e}^0\partial_{\boldsymbol n}\zeta_1 {\rm d}s \nonumber\\
  & \le \sum\limits_{e\subset \partial K} \|\Delta \zeta_1\|_{0,e} \|\Pi_{k-2,e}^0\partial_{\boldsymbol n}\zeta_1\|_{0,e}.\label{zeta1}
\end{align}
According to the trace inequality and the relations in \eqref{scalingzeta0}, one has
\begin{align*}
\|\Pi_{k-2,e}^0\partial_{\boldsymbol n}\zeta_1\|_{0,e}
  & \le \|\Pi_{k-2,e}^0\partial_{\boldsymbol n}\zeta_0\|_{0,e} + \|\Pi_{k-2,e}^0\partial_{\boldsymbol n}v\|_{0,e}   \\
  & \lesssim h_K^{1/2}|\zeta_0|_{2,K} + h_K^{-1/2}|\zeta_0|_{1,K}+ \|\Pi_{k-2,e}^0\partial_{\boldsymbol n}v\|_{0,e} \\
  & \lesssim h_K^{-1}\|v\|_{0,\partial K} + \|\Pi_{k-2,e}^0\partial_{\boldsymbol n}v\|_{0,e}
    \lesssim h_K^{-3/2}|\!|\!|v|\!|\!|_K.
\end{align*}
Applying Lemma \ref{rigA2} to get
$\|\Delta \zeta_1\|_{0,e} \lesssim h_K^{-1/2}\|\Delta \zeta_1\|_{0,K} \lesssim h_K^{-1/2}|\zeta_1|_{2,K}$.
Substituting the above equations into \eqref{zeta1}, we then obtain
$| \zeta_1 |_{2,K} \lesssim h_K^{ - 2}|\!|\!|v|\!|\!|_K$.
By the Poincar\'{e}-Friedrichs inequality for $H^2$ functions and the trace inequality, one has
\begin{align*}
\|\zeta_1\|_{0,K}
& \lesssim h_K^2|\zeta_1|_{2,K} + \Big|\int_{\partial K} \zeta_1 {\rm d}s \Big|+ h_K\Big|\int_{\partial K} \nabla \zeta_1 {\rm d}s \Big| \\
& \lesssim h_K^2|\zeta_1|_{2,K} + h_K^{3/2}\|\nabla \zeta_1 \|_{0,\partial K}
 \lesssim h_K^2|\zeta_1|_{2,K} + h_K|\zeta_1|_{1,K}.
\end{align*}
Given any $\varepsilon >0$, from the integration by parts and the Young's inequality one has
\[|\zeta_1|_{1,K}^2 = -(\Delta \zeta_1, \zeta_1)_K \le |\zeta_1|_{2,K}\|\zeta_1\|_{0,K}
  \lesssim \varepsilon^2 h_K^{-2}\|\zeta_1\|_{0,K}^2 + C(\varepsilon)h_K^2|\zeta_1|_{2,K}^2,\]
which together with the previous estimate and the absorbing technique yields
\begin{equation}\label{zeta1L2}
  \|\zeta_1\|_{0,K} \lesssim h_K^2 |\zeta_1|_{2,K} \lesssim |\!|\!|v|\!|\!|_K.
\end{equation}

Step 4: Let $\eta = \zeta + \zeta_2$, where $\zeta_2\in H^2(K)$ is a $C^0$-continuous nonconforming virtual element function satisfying
\[
\Delta^2\zeta_2 \in \mathbb{P}_l(K), \quad  \zeta_2|_e = \Pi_{k - 2,e}^0(\partial_{\boldsymbol n}\zeta_2|_e) = 0, \quad \Delta \zeta_2|_e \in \mathbb{P}_{k-2}(e), ~~e\subset \partial K.
\]
This means $\zeta_2$ has the following vanishing d.o.f.s: values at vertices, moments of $\zeta_2$ up to $k-2$ on edges and moments of $\partial_{\boldsymbol n}\zeta_2$ up to $k-2$ on edges. Additional requirement is
$\Pi_l^0\eta  = \Pi_l^0v$ or $\Pi_l^0\zeta_2  = \Pi_l^0(v-\zeta)$,
which determines the remaining d.o.f.s of $\zeta_2$: moments up to $l$ on $K$. One easily finds that
\[\eta|_e = v|_e, \quad \Pi_{k - 2,e}^0(\partial_{\boldsymbol n}\eta|_e) = \Pi_{k - 2,e}^0(\partial_{\boldsymbol n}v|_e),~~e\subset \partial K.\]

Integrating by parts and using the inverse inequality for polynomials and the equations \eqref{zeta0L2} and \eqref{zeta1L2}, we obtain
\begin{align*}
  |\zeta_2|_{2,K}^2
  & = \int_K \Delta^2\zeta_2\zeta_2 {\rm d}x = \int_K \Delta^2\zeta_2 \Pi_l^0\zeta_2 {\rm d}x\lesssim \|\Delta^2\zeta_2\|_{0,K}\|\Pi_l^0\zeta_2\|_{0,K} \\
  & \lesssim h_K^{-2}|\zeta_2|_{2,K}(\|\Pi_l^0v\|_{0,K} + \|\zeta_0\|_{0,K}+\|\zeta_1\|_{0,K}) \lesssim h_K^{-2}|\zeta_2|_{2,K}|\!|\!|v|\!|\!|_K,
\end{align*}
or
$|\zeta_2|_{2,K}\lesssim h_K^{-2}|\!|\!|v|\!|\!|_K$.

Taking $w=v_K: = \zeta_0+\zeta_1+\zeta_2$, one can check that $w\in S(K)$. The proof is completed by combining the estimates of each term of $v_K$ in $H^2$ semi-norm.
\end{proof}

Next we establish the inverse inequalities with respect to the continuous $L^2$ norm.

\begin{theorem} \label{cor:invL2}
There hold the following inverse inequalities
\begin{equation}
| v |_{1,K} \lesssim h_K^{-1}\|v\|_{0,K}, \quad | v |_{2,K} \lesssim h_K^{ - 2}\|v\|_{0,K},\quad v \in V_{k,l}(K). \label{inverseL2plate}
\end{equation}
\end{theorem}
\begin{proof}
It suffices to prove the second inequality in \eqref{inverseL2plate}. The first one follows from this inequality and the interpolation inequality (3.3) in  \cite{Huang-Yu-2021} directly. In view of Lemma \ref{lem:invH2mesh}, we only need to verify $|\!|\!|v|\!|\!|_K \lesssim \| v \|_{0,K}$. It is evident that
$\|\Pi_{k - 4}^0v\|_{0,K} \lesssim \|v\|_{0,K}$.
According to the trace inequality, we have
\begin{equation}\label{invbd1}
h_K^{ - 3/2}\| v \|_{0,\partial K} \lesssim h_K^{ - 1}| v |_{1,K} + h_K^{ - 2}\| v \|_{0,K}~ \lesssim \varepsilon | v |_{2,K} + C(\varepsilon ) h_K^{ - 2}\| v \|_{0,K}.
\end{equation}
Similarly,
\begin{equation}\label{invbd2}
h_K^{ - 1/2}\| \partial_{\boldsymbol n}v \|_{0,\partial K}\lesssim \varepsilon | v |_{2,K} + C(\varepsilon )h_K^{ - 2}\| v \|_{0,K}.
\end{equation}
The proof is completed by absorbing $\varepsilon | v |_{2,K}$.
\end{proof}

We are in a position to establish the norm equivalence.
\begin{theorem}\label{thm:normeqc0}
For all $v\in V_{k,l}(K)$, there hold
\begin{align}
h_K^{-1}\|v\|_{0,K}  & \eqsim \|\boldsymbol\chi(v)\|_{l^2}, \label{C0L2chi} \\
|v - \Pi_k^\Delta v|_{2,K} & \eqsim h_K^{-1}\|\boldsymbol\chi (v - \Pi_k^\Delta v)\|_{l^2}, \label{C0L2H2}
\end{align}
where $\Pi_k^\Delta$ is the elliptic projector associated with $V_{k,l}(K)$.
\end{theorem}

\begin{proof}
The lower bound of the $l^2$-$L^2$ estimate:
 $\| \boldsymbol\chi (v) \|_{l^2}  \lesssim  h_K^{-1}\| v \|_{0,K}$
  is a direct manipulation. For details, one can refer to Lemma 3.6 in  \cite{Huang-Yu-2021}.

The inverse inequality in terms of degrees of freedom
$|v|_{2,K} \lesssim h_K^{-1}\|\boldsymbol\chi (v)\|_{l^2}$
follows from Lemma \ref{lem:invH2mesh} and the relations \eqref{norm12}-\eqref{norm4}.

 We now establish a Poincar\'{e}-Friedrichs inequality in terms of degrees of freedom. Applying the Poincar\'{e}-Friedrichs inequality for $H^2$ functions, one gets
\[h_K^{ - 2}\| v \|_{0,K}
  \lesssim | v |_{2,K} + h_K^{ - 2}\Big| \int_{\partial K} v {\rm d}s \Big| + h_K^{-1}\Big| \int_{\partial K} \nabla v {\rm d}s \Big|.\]
Note that
\begin{align}
  \int_{\partial K} \nabla v {\rm d}s
   = \sum\limits_{e \subset \partial K} {\boldsymbol n}_e\int_e \partial_{\boldsymbol n}v {\rm d}s  + \sum\limits_{i = 1}^{N_v} {\boldsymbol t}_{e_i}(v(z_{i + 1}) - v(z_i)), \label{compugrad}
\end{align}
and $k\ge 2$. We then obtain from the definition of d.o.f.s that
$h_K^{ - 2}\| v \|_{0,K} \lesssim | v |_{2,K} + h_K^{-1}\|\boldsymbol\chi (v) \|_{l^2}$,
which together with the inverse inequality established in Step 2 yields the upper bound of \eqref{C0L2chi}.

The lower bound of \eqref{C0L2H2} can be derived by combining the inverse inequality \eqref{inverseL2plate}, Poincar\'{e}-Friedrichs inequality and \eqref{C0L2chi}.
\end{proof}

\begin{theorem}\label{thm:interpC0}
For every $K\in \mathcal{T}_h$ and for every $v\in H^s(K)$ with $2\le s \le k+1$, there holds
$|v - I_Kv|_{m,K} \lesssim h_K^{s-m} | v |_{s,K}$ for $m = 0,1,2,$
where $I_K$ is the nodal interpolation operator mapping into $V_{k,l}$.
\end{theorem}
\begin{proof}
Using the norm equivalence in Theorem \ref{thm:normeqc0}, we have
\begin{equation}\label{Ikvstability}
 \| I_Kv \|_{0,K} \eqsim h_K\| \boldsymbol{\chi }(I_Kv) \|_{l^2} = h_K\| \boldsymbol{\chi }(v) \|_{l^2}.
 \end{equation}
A direct manipulation results in the following stability estimate
$\| I_Kv \|_{0,K} \lesssim \| v \|_{0,K} + h_K| v |_{1,K} + h_K^2|v|_{2,K}$
for the interpolation operator. The desired estimates follow from the Bramble-Hilbert argument.
\end{proof}

\section{Korn's inequality and norm equivalence} \label{sect:Korn-Norm}

\subsection{Korn's inequality}
We will establish a Korn-type inequality in Theorem \ref{thm:KornVh}. To figure out the geometric dependence of the hidden constants, we first recall a classic Korn's inequality with respect to star-shaped domains.

\begin{lemma}\label{KornRot}
   Let $D\subset \mathbb{R}^2$ be a bounded domain of diameter $h$, which is star-shaped with respect to a disc of radius $\rho$. Then for any $\boldsymbol v\in \boldsymbol H^1(D)$ satisfying $\int_D{\nabla  \times \boldsymbol v} {\rm d}x=0$, there holds
  $| \boldsymbol v |_{1,D} \lesssim \| \boldsymbol\varepsilon (\boldsymbol v) \|_{0,D}$,
  where the hidden constant only depends on the aspect ratio $h/\rho$.
\end{lemma}
\begin{proof}
According to Theorem 2.3 in  \cite{Costabel-Dauge-2015}, we have the following gradient estimate
\begin{equation}\label{htauv}
\| w \|_{0,D} \le \beta{(D)}^{-1}\| \nabla w \|_{ - 1,D},\quad w \in L_0^2(D)= \Big\{ {v \in L^2(D):\int_D v {\rm d}x = 0} \Big\},
\end{equation}
where  $\beta (D) \ge \rho/(2h)$.
Given $\boldsymbol v \in \boldsymbol{H}^1(D)$ with $\int_D \nabla  \times \boldsymbol v {\rm d}x=0$, taking $w = \nabla  \times \boldsymbol v = \partial_1v_2 - \partial_2v_1$ in \eqref{htauv}, we obtain
$\| \nabla  \times \boldsymbol v \|_{0,D} \le \frac{2h}{\rho} \| \nabla (\nabla  \times \boldsymbol v) \|_{ - 1,D}$.
Observing the equality \eqref{combieps}, one has
$\| \nabla  \times \boldsymbol v \|_{0,D} \lesssim \| \nabla \boldsymbol\varepsilon (\boldsymbol v) \|_{ - 1,D} \lesssim \| \boldsymbol\varepsilon (\boldsymbol v) \|_{0,D}$.
Therefore,
\[\| \nabla \boldsymbol v \|_{0,D}^2 = \| \boldsymbol\varepsilon (\boldsymbol v) \|_{0,D}^2 + \frac{1}{2}\| \nabla  \times \boldsymbol v \|_{0,D}^2 \lesssim \| \boldsymbol\varepsilon (\boldsymbol v) \|_{0,D}^2.\]
This completes the proof.
\end{proof}

For the star-shaped domain in Lemma \ref{KornRot}, as done in  \cite{Brenner-2003} we define on $D$ an interpolation operator as
$\Pi_D^{RM}: \boldsymbol{H}^1(D) \to {\bf RM}(D),~~\boldsymbol v \mapsto \Pi_D^{RM}\boldsymbol v$,
  equipped with the following constraints
\[
 \Big| \int_D {(\boldsymbol v - \Pi_D^{RM}\boldsymbol v)} {\rm d}x \Big| = 0,\qquad \Big| \int_D {\nabla  \times (\boldsymbol v - \Pi_D^{RM}\boldsymbol v)} {\rm d}x \Big| = 0,\quad \boldsymbol v \in \boldsymbol{H}^1(D).
\]
  Applying the Korn's inequality in Lemma \ref{KornRot} we obtain
  \begin{equation}\label{PRKb}
| \boldsymbol v - \Pi_D^{RM}\boldsymbol v |_{1,D} \lesssim \| \boldsymbol\varepsilon (\boldsymbol v - \Pi_D^{RM}\boldsymbol v) \|_{0,D} = \| \boldsymbol\varepsilon(\boldsymbol v) \|_{0,D},
\end{equation}
which along with the Poincar\'{e}-Friedrichs inequality yields
\begin{equation}\label{PRKc}
\| \boldsymbol v - \Pi_D^{RM}\boldsymbol v \|_{0,D} \lesssim h_D| \boldsymbol v - \Pi_D^{RM}\boldsymbol v |_{1,D} \lesssim h_D\| \boldsymbol\varepsilon (\boldsymbol v) \|_{0,D},
\end{equation}
where $h_D$ is the diameter of $D$ and the hidden constant only depends on the aspect ratio $h/\rho$.

It is easy to check that $(\cdot,\cdot)_D$ and $(\cdot,\cdot)_{\partial D}$ are inner products on ${\bf RM}(D)$, where $D$ is a bounded domain in $\mathbb{R}^2$. In the following, we denote by $Q_D^{RM}$ and $Q_{\partial D}^{RM}$ the $L^2$ projection operators from $\boldsymbol{H}^1(D)$ onto ${\bf RM}(D)$ in terms of the inner products $(\cdot,\cdot)_D$ and $(\cdot,\cdot)_{\partial D}$, respectively.

\begin{lemma}
 Suppose $K$ is a given element in $\mathcal{T}_h$ and $\mathcal{T}_K$ is the associated virtual triangulation. Then for any $\tau \in \mathcal{T}_K$, there hold the following estimates
 \begin{align*}
 & | \boldsymbol v - Q_\tau ^{RM}\boldsymbol v |_{1,\tau } + h_\tau ^{ - 1}\| \boldsymbol v - Q_\tau ^{RM}\boldsymbol v \|_{0,\tau } \lesssim \| \boldsymbol\varepsilon (\boldsymbol v) \|_{0,\tau },\quad \boldsymbol v \in \boldsymbol{H}^1(\tau ),\\
 & | \boldsymbol v - Q_{\partial \tau }^{RM}\boldsymbol v |_{1,\tau } + h_\tau ^{ - 1}\| \boldsymbol v - Q_{\partial \tau }^{RM}\boldsymbol v \|_{0,\tau } \lesssim \| \boldsymbol\varepsilon (\boldsymbol v) \|_{0,\tau },\quad \boldsymbol v \in \boldsymbol{H}^1(\tau ).
 \end{align*}
\end{lemma}
\begin{proof}
  Step 1: According to the minimization property of $L^2$ projections and the estimate \eqref{PRKc}, we immediately obtain
\begin{equation}\label{L2tau}
\| \boldsymbol v - Q_\tau ^{RM}\boldsymbol v \|_{0,\tau } \le \| \boldsymbol v - \Pi _\tau ^{RM}\boldsymbol v \|_{0,\tau } \lesssim h_\tau \| \boldsymbol\varepsilon (\boldsymbol v) \|_{0,\tau }.
\end{equation}
Applying the inverse inequality for polynomials yields
\[\begin{aligned}
  | \boldsymbol v - Q_\tau ^{RM}\boldsymbol v |_{1,\tau }
  &\le | \boldsymbol v - \Pi _\tau ^{RM}\boldsymbol v |_{1,\tau } + | \Pi _\tau ^{RM}\boldsymbol v - Q_\tau ^{RM}\boldsymbol v |_{1,\tau } \hfill \\
  & = | \boldsymbol v - \Pi _\tau ^{RM}\boldsymbol v |_{1,\tau } + | Q_\tau ^{RM}(\boldsymbol v - \Pi _\tau ^{RM}\boldsymbol v) |_{1,\tau } \hfill \\
  & \lesssim | \boldsymbol v - \Pi _\tau ^{RM}\boldsymbol v |_{1,\tau } + h_\tau ^{ - 1}\| \boldsymbol v - Q_\tau ^{RM}\boldsymbol v \|_{0,\tau }. \hfill \\
\end{aligned} \]
The first inequality follows from \eqref{PRKb} and \eqref{L2tau}.

Step 2: In a similar manner, we have
$\| \boldsymbol v - Q_{\partial \tau }^{RM}\boldsymbol v \|_{0,\partial \tau } \le \| \boldsymbol v - \Pi _\tau ^{RM}\boldsymbol v \|_{0,\partial \tau }$.
The conjunction of the trace inequality, \eqref{PRKb} and \eqref{PRKc} implies that
\begin{align*}
  \| \boldsymbol v - Q_{\partial \tau }^{RM}\boldsymbol v \|_{0,\partial \tau }
    \lesssim h_\tau ^{ - 1/2}\| \boldsymbol v - \Pi _\tau ^{RM}\boldsymbol v \|_{0,\tau } + h_\tau ^{1/2}| \boldsymbol v - \Pi _\tau ^{RM}\boldsymbol v |_{1,\tau }
   \lesssim h_\tau ^{1/2}\| \boldsymbol\varepsilon (\boldsymbol v) \|_{0,\tau }.
\end{align*}
Observing the fact that $\| \cdot \|_{\tau }$ and $\| \cdot \|_{\partial \tau }$ are norms on ${\bf RM}(\tau )$, we have, by the scaling arguments and the equivalence of norms on finite-dimensional vector spaces, that
\begin{align*}
  \| \boldsymbol v - Q_{\partial \tau }^{RM}\boldsymbol v \|_{0,\tau }
  &\le \| \boldsymbol v - \Pi _\tau ^{RM}\boldsymbol v \|_{0,\tau } + \| Q_{\partial \tau }^{RM}(\boldsymbol v - \Pi _\tau ^{RM}\boldsymbol v) \|_{0,\tau }  \\
  & \le \| \boldsymbol v - \Pi _\tau ^{RM}\boldsymbol v \|_{0,\tau } + h_\tau ^{1/2}\| Q_{\partial \tau }^{RM}(\boldsymbol v - \Pi _\tau ^{RM}\boldsymbol v) \|_{0,\partial \tau }  \\
  & \lesssim \| \boldsymbol v - \Pi _\tau ^{RM}\boldsymbol v \|_{0,\tau } + h_\tau ^{1/2}\| \boldsymbol v - \Pi _\tau ^{RM}\boldsymbol v \|_{0,\partial \tau }
   \lesssim h_\tau \| \boldsymbol\varepsilon (\boldsymbol v) \|_{0,\tau }.
\end{align*}
Applying the inverse inequality for polynomials, we obtain
\[\begin{aligned}
  | \boldsymbol v - Q_{\partial \tau }^{RM}\boldsymbol v |_{1,\tau }
  &\le | \boldsymbol v - \Pi _\tau ^{RM}\boldsymbol v |_{1,\tau } + | Q_{\partial \tau }^{RM}(\boldsymbol v - \Pi _\tau ^{RM}\boldsymbol v) |_{1,\tau } \hfill \\
  & \lesssim | \boldsymbol v - \Pi _\tau ^{RM}\boldsymbol v |_{1,\tau } + h_\tau ^{ - 1}\| \boldsymbol v - Q_{\partial \tau }^{RM}\boldsymbol v \|_{0,\tau }, \hfill \\
\end{aligned} \]
which implies the second inequality.
\end{proof}

We remark that the projections $ Q_\tau ^{RM}\boldsymbol v$ and  $Q_{\partial \tau} ^{RM}\boldsymbol v$ can be naturally extended to the element $K$. In view of the technique of proving the Bramble-Hilbert estimates under relaxed shape regularity assumptions in \cite{Brenner-Guan-Sung-2017} (see Lemma 4.1 there), we further have
\begin{align}
& | \boldsymbol v - Q_\tau ^{RM}\boldsymbol v |_{1,K } + h_K ^{ - 1}\| \boldsymbol v - Q_\tau ^{RM}\boldsymbol v \|_{0,K } \lesssim \| \boldsymbol\varepsilon (\boldsymbol v) \|_{0,K },\quad \boldsymbol v \in \boldsymbol{H}^1(K ), \label{vQRM1} \\
& | \boldsymbol v - Q_{\partial \tau }^{RM}\boldsymbol v |_{1,K} + h_K ^{ - 1}\| \boldsymbol v - Q_{\partial \tau }^{RM}\boldsymbol v \|_{0,K} \lesssim \| \boldsymbol\varepsilon (\boldsymbol v) \|_{0,K},\quad \boldsymbol v \in \boldsymbol{H}^1(K).
 \end{align}

\begin{theorem}\label{thm:KornQK}
 Under the assumption {\bf A1}, for any $K\in \mathcal{T}_h$ there hold
 \begin{align*}
 & | \boldsymbol v - Q_K^{RM}\boldsymbol v |_{1,K}+ h_K^{ - 1}\| \boldsymbol v - Q_K^{RM}\boldsymbol v \|_{0,K} \lesssim \| \boldsymbol\varepsilon (\boldsymbol v) \|_{0,K},\quad \boldsymbol v \in \boldsymbol{H}^1(K),\\
 & | \boldsymbol v - Q_{\partial K}^{RM}\boldsymbol v |_{1,K} + h_K^{ - 1}\| \boldsymbol v - Q_{\partial K}^{RM}\boldsymbol v \|_{0,K} \lesssim \| \boldsymbol\varepsilon (\boldsymbol v) \|_{0,K},\quad \boldsymbol v \in \boldsymbol{H}^1(K).
 \end{align*}
\end{theorem}
\begin{proof}
  Noting that
 $\int_K (\boldsymbol v - Q_K^{RM}\boldsymbol v) \mathrm{d}x = \boldsymbol 0$ and $\int_{\partial K} (\boldsymbol v - Q_{\partial K}^{RM}\boldsymbol v) \mathrm{d}s = \boldsymbol 0$, we obtain from the Poincar\'{e}-Friedrichs inequality that
  \[h_K^{ - 1}\left( \| \boldsymbol v - Q_K^{RM}\boldsymbol v \|_{0,K} + \| \boldsymbol v - Q_{\partial K}^{RM}\boldsymbol v \|_{0,K} \right) \lesssim | \boldsymbol v - Q_K^{RM}\boldsymbol v |_{1,K} + | \boldsymbol v - Q_{\partial K}^{RM}\boldsymbol v |_{1,K}.\]
  It reduces to prove
\begin{equation}\label{vQK12}
  | \boldsymbol v - Q_K^{RM}\boldsymbol v |_{1,K} + | \boldsymbol v - Q_{\partial K}^{RM}\boldsymbol v |_{1,K} \lesssim \| \boldsymbol\varepsilon (\boldsymbol v) \|_{0,K},\quad \boldsymbol v \in \boldsymbol{H}^1(K).
\end{equation}
  Referring to \eqref{vQRM1} and the inverse estimates for polynomials, we obtain
  \begin{align}
  | \boldsymbol v - Q_K^{RM}\boldsymbol v |_{1,K}
  &\le | \boldsymbol v - Q_\tau ^{RM}\boldsymbol v |_{1,K} + | Q_\tau ^{RM}\boldsymbol v - Q_K^{RM}\boldsymbol v |_{1,K} \nonumber \\
  & \lesssim \| \boldsymbol\varepsilon (\boldsymbol v) \|_{0,K} + h_K^{ - 1}\| Q_\tau ^{RM}\boldsymbol v - Q_K^{RM}\boldsymbol v \|_{0,K} \nonumber \\
  & \le \| \boldsymbol\varepsilon (\boldsymbol v) \|_{0,K} + h_K^{ - 1}( \| \boldsymbol v - Q_\tau ^{RM}\boldsymbol v \|_{0,K} + \| \boldsymbol v - Q_K^{RM}\boldsymbol v \|_{0,K} ) \nonumber \\
  & \le \| \boldsymbol\varepsilon (\boldsymbol v) \|_{0,K} + h_K^{ - 1}\| \boldsymbol v - Q_\tau ^{RM}\boldsymbol v \|_{0,K} \lesssim \| \boldsymbol\varepsilon (\boldsymbol v) \|_{0,K}, \label{firstQK}
\end{align}
  where in the last step the following minimization property of orthogonal projections is used:
  $\| \boldsymbol v - Q_K^{RM}\boldsymbol v \|_{0,K} \le \| \boldsymbol v - Q_\tau ^{RM}\boldsymbol v \|_{0,K}$.

  From the scaling argument and the equivalence of norms on finite-dimensional vector spaces, one easily gets
  $\| \boldsymbol v \|_{0,K} \lesssim \| \boldsymbol v \|_{0,\partial K}, \boldsymbol v \in {\bf RM}(K)$,
  which along with \eqref{firstQK} and the inverse inequality for polynomials yields
  \begin{align*}
  | \boldsymbol v - Q_{\partial K}^{RM}\boldsymbol v |_{1,K}
  &\le | \boldsymbol v - Q_K^{RM}\boldsymbol v |_{1,K} + | Q_K^{RM}\boldsymbol v - Q_{\partial K}^{RM}\boldsymbol v |_{1,K}  \\
  & \lesssim \| \boldsymbol\varepsilon (\boldsymbol v) \|_{0,K} + h_K^{ - 1}\| Q_{\partial K}^{RM}(\boldsymbol v - Q_K^{RM}\boldsymbol v)\|_{0,K} \\
  & \le \| \boldsymbol\varepsilon (\boldsymbol v) \|_{0,K} + h_K^{ - 1}h_K^{1/2}| Q_{\partial K}^{RM}(\boldsymbol v - Q_K^{RM}\boldsymbol v) \|_{0,\partial K}  \\
  &\lesssim \| \boldsymbol\varepsilon (\boldsymbol v) \|_{0,K} + h_K^{ - 1/2}\| \boldsymbol v - Q_K^{RM}\boldsymbol v \|_{0,\partial K}.
\end{align*}
According to the trace inequality and the Poincar\'{e}-Friedrichs inequality, we have
\begin{align*}
  & | \boldsymbol v - Q_{\partial K}^{RM}\boldsymbol v |_{1,K} \\
  & \lesssim \| \boldsymbol\varepsilon (\boldsymbol v) \|_{0,K} + h_K^{ - 1/2}( h_K^{ - 1/2}\| \boldsymbol v - Q_K^{RM}\boldsymbol v \|_{0,K} + h_K^{1/2}| \boldsymbol v - Q_K^{RM}\boldsymbol v |_{1,K} )  \\
  & \lesssim \| \boldsymbol\varepsilon (\boldsymbol v) \|_{0,K} + | \boldsymbol v - Q_K^{RM}\boldsymbol v |_{1,K}.
\end{align*}
This completes the proof by using \eqref{firstQK}.
\end{proof}

We now present the main result.
\begin{theorem} \label{thm:KornVh}
  For any $\boldsymbol{v} \in \boldsymbol{V}_h(K)$, there holds
  \[| \boldsymbol{\tilde{v}} |_{1,K} \eqsim \| \boldsymbol{\varepsilon} (\boldsymbol{\tilde{v}}) \|_{0,K}, \quad \boldsymbol {\tilde{v}} = \boldsymbol{v} - \Pi_K^1\boldsymbol{v}.\]
\end{theorem}
\begin{proof}
In view of the constraint \eqref{constrainPiC}, we have $Q_{\partial K}^{RM}\Pi_K^1\boldsymbol v = Q_{\partial K}^{RM}\boldsymbol v$ or $Q_{\partial K}^{RM} \boldsymbol {\tilde v} = \boldsymbol 0$ for any $\boldsymbol v\in \boldsymbol V_h(K)$. According to the Korn's inequality in Theorem \ref{thm:KornQK}, we immediately obtain
$  | \boldsymbol{\tilde v} |_{1,K}
 = | \boldsymbol{\tilde v} - Q_{\partial K}^{RM}\boldsymbol{\tilde v} |_{1,K} \lesssim  \| \boldsymbol\varepsilon (\boldsymbol{\tilde v}) \|_{0,K}. $
This completes the proof.
\end{proof}

\subsection{Norm equivalence}

According to Theorem \ref{cor:invL2}, we can obtain the following inverse inequalities for $C^0$-continuous nonconforming virtual element spaces.
\begin{lemma}\label{lem:inverseSG}
For all $\boldsymbol{v} \in \boldsymbol{V}_h(K)$, there holds
$| \boldsymbol{v} |_{i,K} \lesssim h_K^{-i} \| \boldsymbol{v} \|_{0,K}$ for $ i = 1,2$.
\end{lemma}

In view of the constraints of each elliptic projector, we obtain the following Poincar\'{e}-Friedrichs inequalities by using the standard Poincar\'{e}-Friedrichs inequalities for $H^1$ and $H^2$ functions.
\begin{lemma}\label{lem:PFSG}
  There hold the Poincar\'{e}-Friedrichs inequalities
  \begin{align*}
  & \| \boldsymbol v - \Pi_K^1\boldsymbol v \|_{0,K} \lesssim h_K| \boldsymbol v - \Pi_K^1\boldsymbol v |_{1,K},\quad \boldsymbol v \in \boldsymbol V_h(K)\\
  & \| \boldsymbol v - \Pi_K^2\boldsymbol v \|_{0,K} \lesssim h_K^2| \boldsymbol v - \Pi_K^2\boldsymbol v |_{2,K},\quad \boldsymbol v \in \boldsymbol V_h(K).
  \end{align*}
\end{lemma}

Now we are able to develop the norm equivalence in $\boldsymbol V_h(K)$.
\begin{lemma}\label{lem:VhL2}
  For any $\boldsymbol v \in \boldsymbol V_h(K)$, the following norm equivalence formulas hold
  \begin{align}
    h_K^{ - 1}\| \boldsymbol v \|_{0,K}  & \eqsim \|\boldsymbol\chi(\boldsymbol v)\|_{l^2}, \nonumber
    \\
     | \boldsymbol v - \Pi_K^1\boldsymbol v |_{1,K} & \eqsim \|\boldsymbol\chi (\boldsymbol v - \Pi_K^1\boldsymbol v)\|_{l^2}, \label{VhL2H1}
     \\
    | \boldsymbol v - \Pi_K^2\boldsymbol v |_{2,K} & \eqsim h_K^{ - 1}\|\boldsymbol\chi (\boldsymbol v - \Pi_K^2\boldsymbol v)\|_{l^2}, \label{VhL2H2}
  \end{align}
  where $\boldsymbol \chi$ is the d.o.f. vector associated with $\boldsymbol v \in \boldsymbol V_h(K)$.
\end{lemma}
\begin{proof}
Let $i=1,2$. The first relation is straightforward in view of Theorem \ref{thm:normeqc0}. For the later two formulas, we obtain from the first relation that
$\|\boldsymbol\chi (\boldsymbol{v} - \Pi_K^i\boldsymbol{v})\|_{l^2} \eqsim h_K^{-1}\| \boldsymbol{v} - \Pi_K^i\boldsymbol{v} \|_{0,K}$.
According to Lemma \ref{lem:PFSG}, there holds
$\| \boldsymbol{v} - \Pi_K^i\boldsymbol{v} \|_{0,K} \lesssim h_K^i | \boldsymbol{v} - \Pi_K^i\boldsymbol{v} |_{i,K}$.
Using Lemma \ref{lem:inverseSG} we further get
$| \boldsymbol{v} - \Pi_K^i\boldsymbol{v} |_{i,K} \lesssim h_K^{-i} \| \boldsymbol{v} - \Pi_K^i\boldsymbol{v} \|_{0,K}$.
The above equations lead to \eqref{VhL2H1} and \eqref{VhL2H2} readily.
\end{proof}

Based on Lemma \ref{lem:VhL2}, we are in a position to verify the following reformulated norm equivalence formulas \eqref{mod1} and \eqref{mod2} for our virtual element methods of the strain gradient elasticity.
\begin{theorem}
  For all $\boldsymbol v \in \boldsymbol V_h(K)$, there hold
\begin{align*}
  \| \boldsymbol\varepsilon (\boldsymbol v - \Pi_K^1\boldsymbol v) \|_{0,K} & \eqsim \| \boldsymbol\chi (\boldsymbol v - \Pi_K^1\boldsymbol v) \|_{l^2}, \\
  \| \nabla \boldsymbol\varepsilon (\boldsymbol v - \Pi_K^2\boldsymbol v) \|_{0,K} & \eqsim h_K^{-1}\| \boldsymbol\chi (\boldsymbol v - \Pi_K^2\boldsymbol v) \|_{l^2}.
\end{align*}
\end{theorem}
\begin{proof}
According to Theorem \ref{thm:KornVh}, we have
$\| \boldsymbol\varepsilon (\boldsymbol v - \Pi_K^1\boldsymbol v) \|_{0,K} \eqsim | \boldsymbol v - \Pi_K^1 \boldsymbol v |_{1,K}.$
The first relation follows from \eqref{VhL2H1} of Lemma \ref{lem:VhL2}. Observing the equality \eqref{combieps}, one immediately obtains
$\| \nabla \boldsymbol\varepsilon (\boldsymbol v - \Pi_K^2\boldsymbol v) \|_{0,K} \eqsim | \boldsymbol v - \Pi_K^2 \boldsymbol v |_{2,K}$,
which completes the proof by using \eqref{VhL2H2} in Lemma \ref{lem:VhL2}.
\end{proof}

\section{Error analysis}

\subsection{An abstract Strang-type lemma}

To derive the robustness with respect to the Lam\'{e} constant $\lambda$, we first establish the crucial commutative relations in the virtual element spaces for locking-free analysis as in elastic problems in  \cite{Beirao-Brezzi-Marini-2013}.

\begin{lemma}\label{lem:divProj}
Let $\boldsymbol{u}\in \boldsymbol{H}^2(\Omega)$ and let $I_h\boldsymbol{u}$ be the interpolation of $\boldsymbol{u}$ in the lifting space ${\boldsymbol{V}}_h$. Then we have
\begin{enumerate}
  \item[(1)] There exists $\boldsymbol{u}_{\mathcal{I}} \in {\boldsymbol{V}}_h$ such that
\begin{align}
 \Pi_{k - 1}^0{\rm div}\boldsymbol{u}_{\mathcal{I}} = \Pi_{k - 1}^0{\rm div}\boldsymbol{u}, \quad
 \Pi_{k-2}^0\nabla {\rm div}\boldsymbol{u}_{\mathcal{I}} = \Pi_{k-2}^0\nabla {\rm div}\boldsymbol{u}. \label{commuteEst}
\end{align}
  \item[(2)] There holds
$|\boldsymbol{u}-\boldsymbol{u}_{\mathcal{I}}|_{\ell,K} \lesssim |\boldsymbol{u}-I_h\boldsymbol{u} |_{\ell,K}$, $\ell = 1,2.$
\end{enumerate}
\end{lemma}
\begin{proof}
Let $\boldsymbol{u}_{\mathcal{I}} = I_h\boldsymbol{u} + \widetilde{\boldsymbol{u}}_h$, where $\widetilde{\boldsymbol{u}}_h\in {\boldsymbol{V}}_h$ is a virtual element function to be determined.

Step 1: We first resolve $\widetilde{\boldsymbol{u}}_h\in {\boldsymbol{V}}_h$, or equivalently define its degrees of freedom.
We choose $\widetilde{\boldsymbol{u}}_h$ with vanishing boundary d.o.f.s, including the values at vertices, moments of $\widetilde{\boldsymbol{u}}_h$ on $e$ and moments of $\partial_{\boldsymbol{n}}\widetilde{\boldsymbol{u}}_h$ on $e$. By the $C^0$-continuity, $\widetilde{\boldsymbol{u}}_h|_{\partial{K}}=\boldsymbol{0}$.
In this case the first equation in \eqref{commuteEst} is equivalent to
\[({\rm div}\widetilde{\boldsymbol{u}}_h, q)_K = ({\rm div}(\boldsymbol{u}-I_h\boldsymbol{u}), q)_K,~~q \in \mathbb{P}_{k-1}(K).\]
Integrating by parts yields
\[({\rm div}\widetilde{\boldsymbol{u}}_h, q)_K =  - ( \widetilde{\boldsymbol{u}}_h, \nabla q )_K + (\widetilde{\boldsymbol{u}}_h \cdot \boldsymbol n, q)_{\partial K} = - ( \widetilde{\boldsymbol{u}}_h, \nabla q )_K.\]
In this way, the equation is reformulated as
\begin{equation}\label{dofKcommute}
- ( \widetilde{\boldsymbol{u}}_h, \nabla q )_K = ({\rm div}(\boldsymbol{u}-I_h\boldsymbol{u}), q)_K,~~q \in \mathbb{P}_{k-1}(K),
\end{equation}
which provides the interior moments of $\widetilde{\boldsymbol{u}}_h$.

Step 2: We proceed to prove the second equation of \eqref{commuteEst}. According to the definition of the d.o.f.s of $\widetilde{\boldsymbol{u}}_h$, it is easy to find that all the boundary terms can be removed in the integration by parts formula in Lemma \ref{lem:d2k}. We then have
\begin{equation}\label{dofnabladiv}
( \nabla {\rm div}\widetilde{\boldsymbol{u}}_h,\boldsymbol{p} )_K
   = \int_K \boldsymbol{Q}_{1i,i}(\boldsymbol{p}) \cdot \widetilde{\boldsymbol{u}}_h {\rm d}x, \quad \boldsymbol{p} \in (\mathbb{P}_{k - 2}(K))^2.
\end{equation}
By integration by parts and using the definition of the interpolation operators, we obtain
\begin{align*}
({\rm div}(\boldsymbol{u}-I_h\boldsymbol{u}), q)_K
 =  - ((\boldsymbol{u}-I_h\boldsymbol{u}), \nabla q )_K + ((\boldsymbol{u}-I_h\boldsymbol{u}) \cdot \boldsymbol n, q)_{\partial K}
 =0, ~~ q \in \mathbb{P}_{k-3}(K).
\end{align*}
From \eqref{dofKcommute} we can get that the interior moments of $\widetilde{\boldsymbol{u}}_h$ up to $k-4$ are zero. Noting that $\boldsymbol{Q}_{1i,i}(\boldsymbol{p})\in (\mathbb{P}_{k-4}(K))^2$, the above result implies that the right-hand side of \eqref{dofnabladiv} is zero, that is, $\Pi_{k-2}^0\nabla {\rm div}\widetilde{\boldsymbol{u}}_h=\boldsymbol{0}$. Hence,
$\Pi_{k-2}^0\nabla {\rm div}\boldsymbol{u}_{\mathcal{I}} = \Pi_{k-2}^0\nabla {\rm div}(I_h\boldsymbol{u})$.
One can check that $\Pi_{k-2}^0\nabla {\rm div}(I_h\boldsymbol{u}) = \Pi_{k-2}^0\nabla {\rm div}\boldsymbol{u}$, which along with the previous equation yields the second one.

Step 3: Let $\ell = 1,2$. By the norm equivalence,
$h_K^{-1}\|\widetilde{\boldsymbol{u}}_h\|_{0,K} \eqsim \|\boldsymbol\chi_4(\widetilde{\boldsymbol{u}}_h)\|_{l^2}$,
where $\boldsymbol\chi_4$ stands for the interior moments. From \eqref{dofKcommute}, the scaling argument gives
\[\|\boldsymbol\chi_4(\widetilde{\boldsymbol{u}}_h)\|_{l^2} \lesssim h_K\|\boldsymbol\chi_4(\Pi_{k-1}^0{\rm div}(\boldsymbol{u}-I_h\boldsymbol{u}))\|_{l^2} \eqsim \|\Pi_{k-1}^0{\rm div}(\boldsymbol{u}-I_h\boldsymbol{u})\|_{0,K}.\]
Hence,
$\|\widetilde{\boldsymbol{u}}_h\|_{0,K} \lesssim h_K |\boldsymbol{u}-I_h\boldsymbol{u}|_{1,K}$.
Applying the inverse inequality to get
$|\widetilde{\boldsymbol{u}}_h|_{\ell,K} \lesssim h_K^{1-\ell} |\boldsymbol{u}-I_h\boldsymbol{u}|_{1,K}$.
According to the Poincar\'{e}-Friedrichs inequality and the definition of the interpolation operators, we obtain
$h_K^{-\ell}\|\widetilde{\boldsymbol{u}}_h\|_{0,K} + |\widetilde{\boldsymbol{u}}_h|_{\ell,K} \lesssim |\boldsymbol{u}-I_h\boldsymbol{u}|_{\ell,K}.$
The triangle inequality leads to
$|\boldsymbol{u}-\boldsymbol{u}_{\mathcal{I}}|_{\ell,K}
 \le  |\boldsymbol{u}-I_h\boldsymbol{u} |_{\ell,K} + |\widetilde{\boldsymbol{u}}_h|_{\ell,K}
 \lesssim |\boldsymbol{u}-I_h\boldsymbol{u}|_{\ell,K}$,
which is the desired estimate.
\end{proof}

We are ready to derive the following Strang-type lemma.
\begin{lemma}\label{lem:StrangLift}
Let $\boldsymbol{u} \in \boldsymbol{H}_0^2(\Omega)$ and $\boldsymbol{u}_h \in \boldsymbol{V}_h$ be the solutions of \eqref{ContinuousP} and \eqref{DiscreteP}, respectively. And let $\boldsymbol{u}_{\mathcal{I}} \in \boldsymbol{V}_h$ be the virtual element function given by Lemma \ref{lem:divProj}. Then for all $\boldsymbol{u}_{\pi}\in (\mathbb{P}_k(\mathcal{T}_h))^2$, there holds
\begin{align}
  \| \boldsymbol{u} - \boldsymbol{u}_h \|_{\iota,h}
  & \lesssim \| \boldsymbol{u} - \boldsymbol{u}_{\mathcal{I}} \|_{\iota,h} + \| \boldsymbol{u} - \boldsymbol{u}_\pi  \|_{\iota,h}  + \| \boldsymbol{f} - \boldsymbol{f}_h \|_{{\boldsymbol{V}}'_h} + E_h \nonumber\\
  &\quad + \lambda ( \| {\rm div}\boldsymbol{u} - \Pi _{k - 1}^0{\rm div}\boldsymbol{u} \|_0 + \iota \| \nabla {\rm div}\boldsymbol{u} - \Pi _{k - 2}^0\nabla {\rm div}\boldsymbol{u} \|_0 )  , \label{StrangtypeLift}
\end{align}
where the hidden constant is independent of $\iota$ and $\lambda$, and
\[\| \boldsymbol{f} - \boldsymbol{f}_h \|_{\boldsymbol{V}'_h} = \mathop {\sup}\limits_{\boldsymbol{v}_h \in {\boldsymbol{V}}_h} \frac{|\langle \boldsymbol{f} - \boldsymbol{f}_h,\boldsymbol{v}_h \rangle|}{\| \boldsymbol{v}_h \|_{\iota,h}},
\quad E_h = \mathop {\sup}\limits_{\boldsymbol{v}_h \in \boldsymbol{V}_h} \frac{|a(\boldsymbol{u},\boldsymbol{v}_h) - (\boldsymbol{f},\boldsymbol{v}_h)|}{\| \boldsymbol{v}_h \|_{\iota,h}}.\]
\end{lemma}
\begin{proof}
Following the similar arguments in \cite{Beirao-Brezzi-Marini-2013,Zhang-Zhao-Yang-2019}, we can deduce that
\begin{align}
  \| \boldsymbol{u} - \boldsymbol{u}_h \|_{\iota,h}
  & \lesssim \| \boldsymbol{u} - \boldsymbol{u}_{\mathcal{I}} \|_{\iota,h} + \| \boldsymbol{u} - \boldsymbol{u}_\pi  \|_{\iota,h}  +  \| \boldsymbol{f} - \boldsymbol{f}_h \|_{\boldsymbol{V}'_h} + E_h\\
  &\quad + \lambda ( \| {\rm div}\boldsymbol{u} - \Pi _{k - 1}^0{\rm div}\boldsymbol{u}_{\mathcal{I}} \|_0 + \iota \| \nabla {\rm div}\boldsymbol{u} - \Pi _{k - 2}^0\nabla {\rm div}\boldsymbol{u}_{\mathcal{I}} \|_0 ) . \label{Strangtype}
\end{align}
The proof is completed by utilizing the commutative relations in Lemma \ref{lem:divProj}.
\end{proof}

\subsection{Error estimates}

The estimate of the load term is given in \cite{Zhao-Chen-Zhang-2016}. For $k\ge 2$, one has
\begin{equation}\label{rhsOpt}
\| \boldsymbol f - \boldsymbol f_h \|_{\boldsymbol V'_h} \lesssim h^k\| \boldsymbol f \|_{k - 1},\quad k \ge 2.
\end{equation}
This ensures the optimal $\mathcal{O}(h^k)$ error bound for data approximation. For the case $k=2$, a variant estimate is
\begin{equation}\label{fhs}
\| \boldsymbol f - \boldsymbol f_h \|_{\boldsymbol V'_h} \lesssim h\| \boldsymbol f \|_0,\quad k = 2,
\end{equation}
which is lower order with respect to $h$ but requires a weaker regularity on $\boldsymbol f$.

We will derive a bound of the consistent term $E_h$ as follows.

\begin{lemma}\label{lem:Eh}
Assume that $\boldsymbol u\in \boldsymbol H^{k+1}(\Omega)$. Then we have the estimate
\[E_h \lesssim \iota h^{k - 1}(| \boldsymbol u |_{k + 1} + \lambda |{\rm div} \boldsymbol u|_k).\]
\end{lemma}
\begin{proof}
 Since $\boldsymbol V_h \subset \boldsymbol H_0^1(\Omega)$, integrating by parts yields
\begin{align*}
  (\boldsymbol{f},\boldsymbol{v}_h)
   = ( - {\rm div}\boldsymbol{\tilde \sigma}(\boldsymbol{u}),\boldsymbol{v}_h) = (\boldsymbol{\tilde \sigma}(\boldsymbol{u}),\boldsymbol\varepsilon (\boldsymbol{v}_h))
   = {\rm I}_1 + \iota^2{\rm I}_2,
\end{align*}
where
${\rm I}_1 = (\boldsymbol{\sigma}(\boldsymbol{u}),\boldsymbol\varepsilon (\boldsymbol{v}_h))$ and ${\rm I}_2 =  - (\Delta \boldsymbol{\sigma}(\boldsymbol{u}),\boldsymbol\varepsilon (\boldsymbol{v}_h))$.
By the definition of the tensors $\boldsymbol{\sigma}$ and $\boldsymbol\varepsilon$, the first term can be written as
${\rm I}_1 = 2\mu (\boldsymbol\varepsilon(\boldsymbol u),\boldsymbol\varepsilon (\boldsymbol v_h)) + \lambda ({\rm div}\boldsymbol u,{\rm div}\boldsymbol v_h)$.
For the second one, we have, by integration by parts, that
\begin{align*}
  {{\rm I}_2}
  & =  - (\Delta {\boldsymbol{\sigma}_{ij}}(\boldsymbol u),\boldsymbol\varepsilon_{ij}(\boldsymbol v_h)) \\
  & = \sum\limits_{K \in \mathcal{T}_h} \int_K {\nabla {\boldsymbol{\sigma}_{ij}}(\boldsymbol u) \cdot \nabla \boldsymbol\varepsilon_{ij}(\boldsymbol v_h)} {\rm d}x  - \sum\limits_{K \in \mathcal{T}_h} \int_{\partial K} \partial_{\boldsymbol n}{\boldsymbol{\sigma}_{ij}}(\boldsymbol u) \boldsymbol\varepsilon_{ij}(\boldsymbol v_h){\rm d}s  \\
  & = (\nabla \boldsymbol{\sigma}(\boldsymbol u),\nabla \boldsymbol\varepsilon (\boldsymbol v_h)) - \sum\limits_{K \in \mathcal{T}_h} \int_{\partial K} \partial_{\boldsymbol n}\boldsymbol{\sigma}_{ij}(\boldsymbol u) \boldsymbol\varepsilon_{ij}(\boldsymbol v_h){\rm d}s .
\end{align*}
It is easy to check that
$(\nabla \boldsymbol{\sigma}(\boldsymbol u),\nabla \boldsymbol\varepsilon (\boldsymbol v_h)) = 2\mu (\nabla \boldsymbol\varepsilon (\boldsymbol u),\nabla \boldsymbol\varepsilon (\boldsymbol v_h)) + \lambda (\nabla {\rm div}\boldsymbol u,\nabla {\rm div}\boldsymbol v_h)$.
Consequently,
\begin{equation}\label{af1}
a(\boldsymbol{u},\boldsymbol{v}_h) - (\boldsymbol{f},\boldsymbol{v}_h) = \iota^2\sum\limits_{K \in \mathcal{T}_h} {\int_{\partial K} \partial_{\boldsymbol n}{{\sigma}_{ij}}(\boldsymbol{u}) \varepsilon_{ij}(\boldsymbol{v}_h){\rm d}s} .
\end{equation}
In view of the symmetry of each tensor and $\partial_iv = n_i\partial_{\boldsymbol n}v + t_i\partial_{\boldsymbol t}v$, we immediately obtain
\begin{align}
  &\sum\limits_{K \in \mathcal{T}_h} \int_{\partial K} \partial_{\boldsymbol n}{\sigma}_{ij}(\boldsymbol{u}) \varepsilon_{ij}(\boldsymbol{v}_h){\rm d}s  \nonumber \\
   =& \sum\limits_{K \in \mathcal{T}_h} \int_{\partial K} {\partial_k}{\sigma}_{ij}(\boldsymbol{u})n_k {\partial_j}v_i{\rm d}s = \sum\limits_{K \in \mathcal{T}_h} {\int_{\partial K} {{\partial_k}{{\sigma}_{ij}}(\boldsymbol{u})n_k} (n_j\partial_{\boldsymbol n}v_i +t_j\partial_{\boldsymbol t}v_i){\rm d}s}  \nonumber \\
   =& \sum\limits_{K \in \mathcal{T}_h} \int_{\partial K} {\partial_k}{\sigma}_{ij}(\boldsymbol{u})n_kn_j\partial_{\boldsymbol n}v_i{\rm d}s = \sum\limits_{e \in \mathcal{E}_h} \int_e n_kn_j{\partial_k}{\sigma}_{ij}(\boldsymbol{u}) [ \partial_{\boldsymbol n}v_i ]{\rm d}s. \label{af2}
\end{align}
Here, $\boldsymbol{v}_h = [v_1,v_2]^\intercal$ and we have used the fact that the contour integration of tangential derivative along the element boundary is zero.

By the weak continuity,
$\int_e [ \partial_{\boldsymbol n}v_i ]p {\rm d}s = 0$ for $p \in \mathbb{P}_{k - 2}(K)$. We then have, by setting $w_i = n_kn_j{\partial_k}{\boldsymbol{\sigma}_{ij}}(\boldsymbol u)$, that
\begin{align}
  \int_e w_i[ \partial_{\boldsymbol n}v_i ] {\rm d}s
 &= \int_e (w_i - \Pi_{k - 2,e}^0w_i)( [ \partial_{\boldsymbol n}v_i
 ] - \Pi_{0,e}^0[ \partial_{\boldsymbol n}v_i ] ) {\rm d}s \nonumber \\
 & \le \|w_i - \Pi_{k - 2,e}^0w_i\|_{0,e}\| [ \partial_{\boldsymbol n}v_i ] - \Pi_{0,e}^0[ \partial_{\boldsymbol n}v_i
 ]\|_{0,e}. \label{af3}
\end{align}
According to the minimization property of $L^2$ projection, we obtain from the trace inequality, for each internal edge $e = \partial K^ +  \cap \partial K^ - $, that
\begin{align}
  \| w_i - \Pi_{k - 2,e}^0w_i \|_{0,e}
  & \le \| w_i - \Pi_{k - 2}^0w_i \|_{0,e} \nonumber \\
  & \lesssim  h^{1/2}| {w_i - \Pi_{k - 2}^0w_i} |_{1,K^ +  \cup K^ -} + h^{ - 1/2}\| w_i - \Pi_{k - 2}^0w_i \|_{0,K^ +  \cup K^ -} \nonumber \\
  & \lesssim h^{k - 3/2}| w_i |_{k - 1,K^ +  \cup K^ -}, \label{af4}
\end{align}
and similarly,
\begin{equation}\label{af5}
\| [ \partial_{\boldsymbol n}v_i ] - \Pi_{0,e}^0[ \partial_{\boldsymbol n}v_i ] \|_{0,e} \lesssim h^{1/2}| \boldsymbol{v}_h |_{2,K^ +  \cup K^ -}.
\end{equation}
Noting that
$w_i = n_kn_j\partial_k\sigma_{ij}(\boldsymbol{u}) = n_kn_j( \lambda \partial_k({\rm div}\boldsymbol{u})\delta_{ij} + 2\mu {\partial_k}{\varepsilon_{ij}}(\boldsymbol{u}) )$,
we have
\begin{equation}\label{af6}
| w_i |_{k - 1,K^ +  \cup K^ -} \lesssim h^{k - 3/2}(| \boldsymbol{u} |_{k + 1,K^ +  \cup K^ -} + \lambda | {\rm div}\boldsymbol{u}|_{k,K^ +  \cup K^ -}).
\end{equation}
For boundary edges, the adjustment is obvious since $K^+ = K^- = K$. The proof is completed by combining \eqref{af1}-\eqref{af6}.
\end{proof}

To sum up the above results, we obtain the error estimate for the VEMs described as follows.

\begin{theorem}\label{thm:errLift}
Given $k\ge 2$ and $\boldsymbol{f} \in \boldsymbol{L}^2(\Omega ) \cap \boldsymbol{H}^{k - 1}(\Omega )$, let $\boldsymbol{u} \in \boldsymbol{H}_0^2(\Omega ) \cap \boldsymbol{H}^{k + 1}(\Omega )$ be the solution of the continuous variational problem \eqref{ContinuousP}, and let $\boldsymbol{u}_h\in \boldsymbol{V}_h={\boldsymbol{V}}_h$ be the VEM solution of \eqref{DiscreteP}. Then there holds
\[\| \boldsymbol{u} - \boldsymbol{u}_h \|_{\iota,h}
  \lesssim  \begin{cases}
  (h^k + \iota h^{k - 1})| \boldsymbol{u}|_{k + 1} + h^k\| \boldsymbol{f} \|_{k - 1}+E_\lambda,  \\
  h^{k - 1}(| \boldsymbol{u} |_k + \iota | \boldsymbol{u} |_{k + 1}) + h^k\| \boldsymbol{f} \|_{k - 1}+E_\lambda,
\end{cases}  \]
where $E_\lambda  = \lambda (h^k + \iota h^{k - 1})| {\rm div}\boldsymbol{u} |_k$ and the hidden constant is independent of $\iota$ and $\lambda$.
\end{theorem}
\begin{proof}
We only need to bound each term of the right-hand side of \eqref{StrangtypeLift} in Lemma \ref{lem:StrangLift}. For $\lambda$-dependent terms, we obtain from the error estimates of $L^2$ projections that
\begin{align*}
\lambda ( \| \mathrm{div}\boldsymbol{u} - \Pi _{k - 1}^0\mathrm{div}\boldsymbol{u} \|_0 + \iota\| \nabla \mathrm{div}\boldsymbol{u} - \Pi _{k - 2}^0\nabla \mathrm{div}\boldsymbol{u} \|_0 )
\lesssim \lambda (h^k + \iota h^{k-1})|{\rm div}\boldsymbol{u}|_k = E_\lambda.
\end{align*}

The first term is determined by \eqref{commuteEst}. In fact, by the interpolation error estimate in Theorem \ref{thm:interpC0} and the Dupont-Scott theory, there exists $\boldsymbol{u}_\pi\in (\mathbb{P}_k(K))^2$ such that
\[\| \boldsymbol{u} - \boldsymbol{u}_{\mathcal{I}} \|_{\iota,h} + \| \boldsymbol{u} - \boldsymbol{u}_\pi  \|_{\iota,h}
\lesssim h^{\ell-1}|\boldsymbol{u}|_\ell + \iota h^{\ell'-2}|\boldsymbol{u}|_{\ell'}.\]
If $\ell = \ell' = k+1$, then
$\| \boldsymbol{u} - \boldsymbol{u}_{\mathcal{I}} \|_{\iota,h} + \| \boldsymbol{u} - \boldsymbol{u}_\pi  \|_{\iota,h}
\lesssim (h^k + \iota h^{k-1})|\boldsymbol{u}|_{k+1}$,
while for $\ell = k$, $\ell'=k+1$,
$\| \boldsymbol{u} - \boldsymbol{u}_{\mathcal{I}} \|_{\iota,h} + \| \boldsymbol{u} - \boldsymbol{u}_\pi  \|_{\iota,h}
\lesssim h^{k-1}(|\boldsymbol{u}|_k + \iota |\boldsymbol{u}|_{k+1})$.
The proof is completed by combining the above equations, \eqref{rhsOpt} and Lemma \ref{lem:Eh}.
\end{proof}

The following regularity estimate is established for a bounded convex polygonal domain $\Omega \subset \mathbb{R}^2$ in \cite{Liao-Ming-Xu-2021}: there exists $C$ independent of $\iota$ and $\lambda$ such that
\begin{equation}\label{regularitynotsharp}
\|\boldsymbol{u}\|_3 + \lambda \|{\rm div} \boldsymbol{u} \|_2 \le C \iota^{-2} \|\boldsymbol{f}\|_{-1},
\end{equation}
which leads to
$E_\lambda  = \lambda (h^k + \iota h^{k - 1})| {\rm div}\boldsymbol{u} |_k \le C_f ( ( {h}/{\iota} )^2 + {h}/{\iota} )$
in the lowest order case $k=2$. We remark that the estimate in \eqref{regularitynotsharp} is not sharp, hence it does not imply the robustness with respect to the microscopic parameter $\iota$.

\subsection{Uniform error estimate in the lowest order case}

To derive the uniform error estimate with respect to the parameters $\iota$ and $\lambda$, we must establish a sharper estimate than the one given in \eqref{regularitynotsharp}. Such an expected regularity in fact has been achieved in Ref.~\cite{CHH2022} under two reasonable assumptions by using the techniques from the mixed formulation.
Let us first recall the assumptions.

\begin{enumerate}
  \item[\textbf{H1.}] Suppose that $\boldsymbol{f} \in \boldsymbol{H}^{-1}(\Omega)$. Let $(\boldsymbol{\phi},p) \in \boldsymbol{H}_0^2(\Omega) \times  ( H_0^1(\Omega) \cap L_0^2(\Omega) )$ be the solution of the following problem
 \begin{equation*}
  \begin{cases}
  \Delta^2 \boldsymbol{\phi} + \nabla \Delta p = \boldsymbol{f} \quad & \text{in}~~ \Omega, \\
  \Delta \text{div} \boldsymbol{\phi} = 0 \quad & \text{in}~~ \Omega, \\
  \boldsymbol{\phi} = \partial_{\boldsymbol{n}} \boldsymbol{\phi} = \boldsymbol{0} \quad & \text{on}~~ \partial \Omega.
  \end{cases}
 \end{equation*}
  Then $\boldsymbol{\phi} \in \boldsymbol{H}_0^3(\Omega)$ and $p \in H^2(\Omega)$, and there holds the estimate
  $\|\boldsymbol{\phi}\|_3 + \|p\|_2 \lesssim \|\boldsymbol{f}\|_{-1}$.
\end{enumerate}

The above assumption has been rigorously proved in \cite{CHH2022} for a convex polygon. It is worth pointing out that $- {\rm div}\boldsymbol{\tilde \sigma }(\boldsymbol{u}) = (\iota^2\Delta  - I)\mathcal{L} \boldsymbol{u}$, where $\mathcal{L} \boldsymbol{u} = \mu \Delta \boldsymbol{u} + (\lambda  + \mu )\nabla ({\rm div} \boldsymbol{u})$ is the usual Lam\'e operator, while $\Delta \mathcal{L}\boldsymbol{u} = \boldsymbol{f}$ can be written in the following mixed form:
\[\Delta^2 \boldsymbol{u} +  \nabla \Delta p = \frac{1}{\mu} \boldsymbol{f}, \qquad p = \frac{\mu+\lambda}{\mu}  \text{div} \boldsymbol{u},\]
which plays the role of the Stokes problem in deriving the locking-free a-priori estimate for the linear elasticity problems (see Eq.~(2.23) in \cite{Brenner-Sung-1992}).

\begin{enumerate}
  \item[\textbf{H2.}] Suppose that $\boldsymbol{f} \in \boldsymbol{H}^{-1}(\Omega)$. Let $\boldsymbol{u} \in \boldsymbol{H}_0^2(\Omega)$ be the solution of the following problem
 \begin{equation}\label{auxmixed}
  \begin{cases}
  \Delta (\mathcal{L} \boldsymbol{u}) = \boldsymbol{f} \quad & \text{in}~~ \Omega, \\
  \boldsymbol{u} = \partial_{\boldsymbol{n}} \boldsymbol{u} = \boldsymbol{0} \quad & \text{on}~~ \partial \Omega,
  \end{cases}
 \end{equation}
 where $\mathcal{L} \boldsymbol{u} = \mu \Delta \boldsymbol{u} + (\lambda  + \mu )\nabla ({\rm div} \boldsymbol{u})$ is the usual Lam\'e operator, with $\lambda \in [0, \Lambda]$ and $\mu \in [\mu_0, \mu_1]$ being the Lam\'e constants. Then $\boldsymbol{u} \in \boldsymbol{H}_0^3(\Omega)$ and it admits the estimate
  $\|\boldsymbol{u}\|_3 \lesssim \|\boldsymbol{f}\|_{-1}$,
  where the hidden constant may depend on $\Lambda$, $\mu_0$ and $\mu_1$.
\end{enumerate}

The second assumption is used to deal with the case $\lambda \le \lambda_0$ for some sufficiently large number $\lambda_0$ as in the linear elasticity problems.

Under the above two assumptions, one may build up the regularity of the SGE problem, with the details of the proof shown in  \cite{CHH2022}. In what follows, let $\boldsymbol u^0 \in \boldsymbol H_0^1(\Omega)$ be the solution of the following linear elasticity problem
\begin{equation}\label{reducedProblem}
   - \rm{div}\boldsymbol\sigma (\boldsymbol u^0) = \boldsymbol f \quad  {\rm in}~~\Omega ,  \quad
  \boldsymbol u^0  = \boldsymbol 0\quad  {\rm on}~~\partial \Omega.
\end{equation}
It is well-known that $\boldsymbol{u}_0$ has the regularity
\begin{equation}\label{u0estimate}
\| \boldsymbol{u}^0 \|_2 + \lambda \| {\rm div}\boldsymbol{u}^0 \|_1 \le C_\Omega \| \boldsymbol{f} \|_0,
\end{equation}
when $\Omega \subset \mathbb{R}^2$ is a bounded convex polygonal domain.

\begin{lemma}[\cite{CHH2022}] \label{lem:regularity}
 Let $\Omega \subset \mathbb{R}^2$ be a bounded convex polygonal domain. Assume that $\boldsymbol f \in \boldsymbol{L}^2(\Omega )$ and $\boldsymbol u \in \boldsymbol{H}_0^2(\Omega )$ is the solution of \eqref{Pro:SGE}. Under the assumptions \textbf{H1} and \textbf{H2}, it holds
 \[|\boldsymbol{u}-\boldsymbol{u}_0|_1 + \iota \|\boldsymbol{u}\|_2 + \iota^2 \|\boldsymbol{u}\|_3 + \lambda \|\rm{div} (\boldsymbol{u}-\boldsymbol{u}_0)\| + \lambda \iota |\rm{div} \boldsymbol{u}|_1 + \lambda \iota^2 \|\rm{div} \boldsymbol{u}\|_2 \lesssim \iota^{1/2} \|\boldsymbol{f}\|_0,\]
 where the hidden constant is independent of $\lambda$ and $\iota$.
\end{lemma}

We are ready to derive the $\iota$- and $\lambda$-independent error estimate.
To this end, we first present a lower order estimate for the interpolation error with the proof omitted.
\begin{lemma}\label{lowInterp}
  For any $v\in H^2(K)$, there holds
 $| v - I_Kv |_{1,K} \lesssim h_K^{1/2}\| v \|_{1,K}^{1/2}\| v \|_{2,K}^{1/2}$.
\end{lemma}

\begin{theorem}\label{thm:robust}
Let $k=2$ be the order of the virtual element space. Under the assumption given in Lemma \ref{lem:regularity}, one has
\[\| \boldsymbol u - \boldsymbol u_h \|_{\iota,h} \lesssim h^{1/2}\|\boldsymbol f\|_0,\qquad
\|\boldsymbol{u}_0 - \boldsymbol{u}_h \|_{\iota,h} \lesssim (h^{1/2} + \iota^{1/2})\|\boldsymbol f\|_0, \]
with the hidden constant independent of $\lambda$ and $\iota$.
\end{theorem}
\begin{proof}
We only need to estimate the right-hand side of \eqref{Strangtype} in Lemma \ref{lem:StrangLift} term by term. In the following, $\boldsymbol{u}_{\mathcal{I}}$ is the function given in Lemma \ref{lem:divProj}.

Step 1: From Lemmas \ref{lem:divProj} and \ref{lem:regularity}, one gets
  \[\begin{aligned}
  \iota | \boldsymbol u - \boldsymbol u_{\mathcal{I}} |_{2,h}
  &= \iota | \boldsymbol u - \boldsymbol u_{\mathcal{I}} |_{2,h}^{1/2}| \boldsymbol u - \boldsymbol u_{\mathcal{I}} |_{2,h}^{1/2} \lesssim \iota h^{1/2}| \boldsymbol u |_2^{1/2}| \boldsymbol u |_3^{1/2} \hfill \\
  & \lesssim \iota h^{1/2}\iota ^{ - 1/4}\| \boldsymbol f \|_0^{1/2}\iota ^{ - 3/4}\| \boldsymbol f \|_0^{1/2} = h^{1/2}\| \boldsymbol f \|_0. \hfill \\
\end{aligned} \]
For $\boldsymbol u^0$ given in Lemma \ref{lem:regularity}, we obtain from the triangle inequality and Lemma \ref{lowInterp} that
  \[\begin{aligned}
  | \boldsymbol u - \boldsymbol u_{\mathcal{I}} |_{1,h}
  &\le |(\boldsymbol u - \boldsymbol u^0) - (\boldsymbol u - \boldsymbol u^0)_I|_{1,h} + |\boldsymbol u^0 - (\boldsymbol u^0)_I|_{1,h} \hfill \\
  & \lesssim h^{1/2}\| \boldsymbol u - \boldsymbol u^0 \|_1^{1/2}\| \boldsymbol u - \boldsymbol u^0 \|_2^{1/2} + h| \boldsymbol u^0 |_2. \hfill \\
\end{aligned} \]
Using Lemma \ref{lem:regularity} and observing the regularity \eqref{u0estimate}, we derive
\[\begin{aligned}
  | \boldsymbol u - \boldsymbol u_{\mathcal{I}} |_{1,h}
  &\lesssim h^{1/2}\iota ^{1/4}\| \boldsymbol f \|_0^{1/2}(\iota ^{ - 1/4}\| \boldsymbol f \|_0^{1/2} + \| \boldsymbol f \|_0^{1/2}) + h\| \boldsymbol f \|_0 \hfill \\
  & \lesssim h^{1/2}(1 + \iota ^{1/4})\| \boldsymbol f \|_0 + h\| \boldsymbol f \|_0 \lesssim h^{1/2}\| \boldsymbol f \|_0, \hfill \\
\end{aligned} \]
where we have used the fact that $0<\iota\le 1$ and without loss of generality we set $0<h < 1$. Combining above estimates yields
$\| \boldsymbol u - \boldsymbol u_{\mathcal{I}} \|_{E,h} \lesssim h^{1/2}\| \boldsymbol f \|_0$.
Proceeding in an analogous fashion, we can derive
$\| \boldsymbol u - \boldsymbol u_\pi  \|_{E,h} \lesssim h^{1/2}\| \boldsymbol f \|_0$.

Step 2: The third and fourth terms can be argued as in the proof of Step 1. Take the third term as an example. We obtain from the triangle inequality and \eqref{u0estimate} that
\[\begin{gathered}
  \lambda \| \mathrm{div}\boldsymbol u - \Pi _1^0\mathrm{div}\boldsymbol u \|_0 \hfill \\
   \le \lambda \| (\mathrm{div}\boldsymbol u - \mathrm{div}\boldsymbol u^0) - \Pi _1^0(\mathrm{div}\boldsymbol u - \mathrm{div}\boldsymbol u^0) \|_0 + \lambda \| \mathrm{div}\boldsymbol u^0 - \Pi _1^0\mathrm{div}\boldsymbol u^0 \|_0 \hfill \\
   \lesssim h^{1/2}\lambda \| \mathrm{div}\boldsymbol u - \mathrm{div}\boldsymbol u^0 \|_0^{1/2}\| \mathrm{div}\boldsymbol u - \mathrm{div}\boldsymbol u^0 \|_1^{1/2} + h\lambda \| \mathrm{div}\boldsymbol u^0 \|_1 \hfill \\
   \lesssim h^{1/2}\lambda \| \mathrm{div}\boldsymbol u - \mathrm{div}\boldsymbol u^0 \|_0^{1/2}\| \mathrm{div}\boldsymbol u - \mathrm{div}\boldsymbol u^0\|_1^{1/2} + h\| \boldsymbol f\|_0. \hfill \\
\end{gathered} \]
Using again the regularity of the SGE model to get
\[\begin{gathered}
  \lambda \| \mathrm{div}\boldsymbol u - \Pi _1^0\mathrm{div}\boldsymbol u \|_0
   \lesssim h^{1/2}\iota ^{1/2}\| \boldsymbol f \|_0^{1/2}(\lambda \| \mathrm{div}\boldsymbol u \|_1 + \lambda \| \mathrm{div}\boldsymbol u^0 \|_1)^{1/2} + h\| \boldsymbol f \|_0 \hfill \\
   \le h^{1/2}(\iota ^{1/2} + \iota )^{1/2}\| \boldsymbol f \|_0 + h\| \boldsymbol f \|_0 \lesssim h^{1/2}\| \boldsymbol f \|_0. \hfill \\
\end{gathered} \]

Step 3: The estimate of the right-side hand is given in \eqref{fhs}. It remains to bound the consistency term. As in the proof of Lemma \ref{lem:Eh},
we have
\[a(\boldsymbol u,\boldsymbol v_h) - (\boldsymbol f,\boldsymbol v_h) = \iota ^2 \sum\limits_{e \in \mathcal{E}_h} \int_e n_kn_j{\partial _k}\boldsymbol\sigma _{ij}(\boldsymbol u) [ \partial _{\boldsymbol n}v_i ]\mathrm{d}s.\]
Setting $w_i = n_kn_j{\partial _k}{\boldsymbol{\sigma }_{ij}}(\boldsymbol u)$, from \eqref{af3} one obtains
\begin{align*}
  \int_e w_i[ \partial _{\boldsymbol n}v_i ] \mathrm{d}s
  & \le \|w_i - \Pi _{0,e}^0w_i\|_{0,e}\| [ \partial _{\boldsymbol n}v_i ] - \Pi _{0,e}^0[ \partial _{\boldsymbol n}v_i ]\|_{0,e}  \\
  & \le \|w_i - \Pi _0^Kw_i\|_{0,e}\| [ \partial _{\boldsymbol n}v_i ] - \Pi _{0,e}^0[ \partial _{\boldsymbol n}v_i ]\|_{0,e},
\end{align*}
where we have used the minimization property of $L^2$ projections.
For each internal edge $e = \partial K^ +  \cap \partial K^ - $, the trace inequality yields
$\| w_i - \Pi _0^Kw_i \|_{0,e} \lesssim \| w_i \|_{0,K}^{1/2}\| w_i \|_{1,K}^{1/2}$.
Noting that
$w_i = n_kn_j\partial _k\boldsymbol\sigma _{ij}(\boldsymbol u) = n_kn_j( \lambda \partial _k(\mathrm{div}\boldsymbol u)\delta _{ij} + 2\mu {\partial _k}{\varepsilon _{ij}}(\boldsymbol u) )$,
we have
\[\| w_i - \Pi _0^Kw_i \|_{0,e} \lesssim ( \lambda \| \mathrm{div}\boldsymbol u \|_{1,K} + \| \boldsymbol u \|_{2,K} )^{1/2}( \lambda \| \mathrm{div}\boldsymbol u \|_{2,K} + \| \boldsymbol u \|_{3,K} )^{1/2}.\]
In view of \eqref{af5}, we then obtain
$E_h \lesssim \iota h^{1/2}( \lambda\| \mathrm{div}\boldsymbol u \|_1 + \| \boldsymbol u \|_2 )^{1/2}( \lambda \| \mathrm{div}\boldsymbol u \|_2 + \| \boldsymbol u \|_3 )^{1/2}$.
Along with the estimates in Lemma \ref{lem:regularity}, it implies
$E_h \lesssim  h^{1/2} \|\boldsymbol f\|_0$.

The second estimate follows from the triangle's inequality, the first estimate and the regularities in \eqref{u0estimate} and Lemma \ref{lem:regularity}.
\end{proof}

\section{Numerical examples} \label{sect:numerical}

In this section, we report the performance of our proposed virtual element method with several examples by testing the accuracy and the robustness with respect to the microscopic parameter $\iota$ and the Lam\'{e} constant $\lambda$.
For simplicity, we only consider the lowest order element ($k = 2$). Unless otherwise specified, the domain $\Omega$ is taken as the unit square $(0,1)^2$.

\begin{example}\label{example1}
We first consider the effect of the microscopic parameter $\iota$. The exact solution $\boldsymbol u = (u_1, u_2)^\intercal$ is
\begin{equation}\label{solu1}
\boldsymbol{u} =
\left[ \begin{array}{*{20}{c}}
({\rm e}^{\cos 2\pi x} - {\rm e})({\rm e}^{\cos 2\pi y} - {\rm e}) \\
(\cos 2\pi x - 1)(\cos 4\pi y - 1)
\end{array} \right]
~~ \mbox{or} ~~
\left[ \begin{array}{*{20}{c}}
\iota \left( {\rm e}^{ - x/\iota } + {\rm e}^{ - y/\iota }\right) - x^2y \\
\iota \left( {\rm e}^{ - x/\iota } + {\rm e}^{ - y/\iota } \right) - xy^2
\end{array} \right],
\end{equation}
where the second one has boundary layer.
\end{example}

Let $\boldsymbol u$ be the exact solution of \eqref{ContinuousP} and $\boldsymbol u_h$ the discrete solution of the underlying VEM \eqref{DiscreteP}. Since the VEM solution is not explicitly known inside the polygonal elements, as in \cite{Beirao-Lovadina-Russo-2017} we will evaluate the errors by comparing the exact solution $\boldsymbol u$ with the elliptic projection $\Pi_K^i\boldsymbol u_h (i=1,2)$. In this way, the discrete relative error in terms of the discrete energy norm is quantified by
\[E_\Pi = ( {\sum\limits_{K \in \mathcal{T}_h} (|\boldsymbol u - \Pi_K^1\boldsymbol u_h|_{1,K}^2 + \iota^2 |\boldsymbol u - \Pi_K^2\boldsymbol u_h|_{2,K}^2)}/({|\boldsymbol u|_{1,\Omega}^2 + \iota^2|\boldsymbol u|_{2,\Omega}^2}) )^{1/2}.\]

To test the accuracy of the proposed method we consider a sequence of meshes, which is a Centroidal Voronoi Tessellation of the unit square in 32, 64, 128, 256 and 512 polygons. These meshes are generated by the MATLAB toolbox - PolyMesher introduced in \cite{Talischi-Paulino-Pereira-2012}.
 The convergence orders of the errors for the first solution against the mesh size $h$ are shown in the top half of Table \ref{tab:exam1poly}.
As observed from Table \ref{tab:exam1poly}, the convergence rate appears to be linear when $\iota$ is large, and the VEM ensures the quadratic convergence as $\iota\to 0$, which is consistent with the theoretical prediction in Theorem \ref{thm:errLift}. Moreover, a stable trend of the errors is observed as $\iota$ decreases to zero. The numerical results for the solution with boundary layer are listed in the bottom half. It is evident that only a half-order rate of convergence is observed for small $\iota$.

\begin{table}[htb]
  \centering
\caption{The convergence rate for Example \ref{example1}($\lambda=\mu=1$))}\label{tab:exam1poly}
  \begin{tabular}{ccccccccccccccccc}
  \toprule[0.2mm]
  $\iota \backslash N$ & 32   &64   &128   &256   &512  & Rate\\
  \midrule[0.3mm]
   1e-0  & 1.1832e+00   & 9.6323e-01   & 6.5934e-01   & 4.2552e-01   & 2.6895e-01 & 1.09\\
   1e-1  & 5.2051e-01   & 4.0272e-01   & 3.3624e-01   & 2.6913e-01   & 1.9874e-01 & 0.67\\
   1e-2  & 3.5613e-01   & 1.7987e-01   & 9.8415e-02   & 5.1840e-02   & 2.8723e-02 & 1.81\\
   1e-3  & 3.5205e-01   & 1.7492e-01   & 9.3162e-02   & 4.6501e-02   & 2.3509e-02 & 1.94\\
   1e-4  & 3.5200e-01   & 1.7487e-01   & 9.3107e-02   & 4.6443e-02   & 2.3450e-02 & 1.95\\
   1e-5  & 3.5200e-01   & 1.7487e-01   & 9.3106e-02   & 4.6442e-02   & 2.3450e-02 & 1.95\\
  \bottomrule[0.2mm]
\end{tabular}
\centering
  \begin{tabular}{cccccccccccccc}
   1e-0  & 1.6613e-01   & 1.1638e-01   & 8.0678e-02   & 5.3877e-02   & 3.6209e-02 & 1.10\\
   1e-1  & 1.6077e-01   & 1.3559e-01   & 1.1002e-01   & 8.9844e-02   & 6.7354e-02 & 0.62\\
   1e-2  & 2.5674e-01   & 2.0325e-01   & 1.6700e-01   & 1.3365e-01   & 1.0540e-01 & 0.63\\
   1e-3  & 2.9559e-01   & 2.4860e-01   & 2.1769e-01   & 1.8364e-01   & 1.5178e-01 & 0.47\\
   1e-4  & 2.9589e-01   & 2.4889e-01   & 2.1801e-01   & 1.8427e-01   & 1.5384e-01 & 0.46\\
   1e-5  & 2.9593e-01   & 2.4893e-01   & 2.1804e-01   & 1.8431e-01   & 1.5387e-01 & 0.46\\
  \bottomrule[0.2mm]
\end{tabular}
\end{table}

\begin{example}\label{example3}
We now investigate the impact of Lam\'{e} constant $\lambda$. The exact solution is given by
 \[\boldsymbol u(x,y) = \left[ \begin{array}{*{20}{c}}
  ( - 1 + \cos 2\pi x)\sin2\pi y \\
   - ( - 1 + \cos 2\pi y)\sin2\pi x
\end{array} \right] + \frac{1}{1 + \lambda}\sin \pi x\sin \pi y\left[ \begin{array}{*{20}{c}}
  1 \\
  1
\end{array} \right].\]
\end{example}

We further introduce a relative error in the discrete maximum norm defined by
\[E_\infty  = {\mathop {\max}\limits_{1 \le i \le N} | \boldsymbol u(x_i,y_i) - \boldsymbol u_h(x_i,y_i) |}/{\mathop {\max}\limits_{1 \le i \le N} | \boldsymbol u(x_i,y_i) |},\]
 where $(x_i,y_i)_{1\le i\le N}$ are all vertices in $\mathcal{T}_h$.
Let $\lambda$ vary from $10$ to $10^5$ on a mesh of 100 elements with $\mu = 1$ and $\iota = 10^{-5}$. The relative errors for different values of $\lambda$ are listed in Table \ref{tab:exam3poly}, from which we may conclude that the VEM is locking-free, which is consistent with the theoretical prediction in Theorem \ref{thm:errLift} since the errors are hardly affected by the choice of $\lambda$ when $\lambda \| {\rm div}\boldsymbol u \|_2$ and $\|\boldsymbol f\|_0$ are uniformly bounded with respect to $\lambda$ as given in Table \ref{tab:exam3poly}.


\begin{table}[htb]
  \centering
  \caption{The relative errors for different $\lambda$ on the fixed mesh for Example \ref{example3} ($\mu = 1$, $\iota = 10^{-5}$)}\label{tab:exam3poly}
  \begin{tabular}{lllllllllllllllllll}
  \toprule[0.2mm]
  $\lambda$ &$10$   &$10^2$   &$10^3$   &$10^4$  &$10^5$\\
  \midrule[0.3mm]
   $\lambda \| {\rm div}\boldsymbol u \|_2$ &4.0910e+01 &4.4556e+01 &4.4956e+01 &4.4997e+01 &4.5001e+01\\[1mm]
   $\|\boldsymbol f \|_0$ &6.9227e+01 &6.9101e+01 &6.9089e+01 &6.9087e+01 &6.9087e+01\\[1mm]
   $E_\Pi$      &9.5502e-02 &9.5367e-02 &9.5359e-02 &9.5358e-02 &9.5358e-02\\[1mm]
   $E_\infty$   &9.2877e-02 &9.2731e-02 &9.2721e-02 &9.2721e-02 &9.2720e-02\\
  \bottomrule[0.2mm]
\end{tabular}
\end{table}

\begin{figure}[htb]
  \centering
 \includegraphics[scale=0.35]{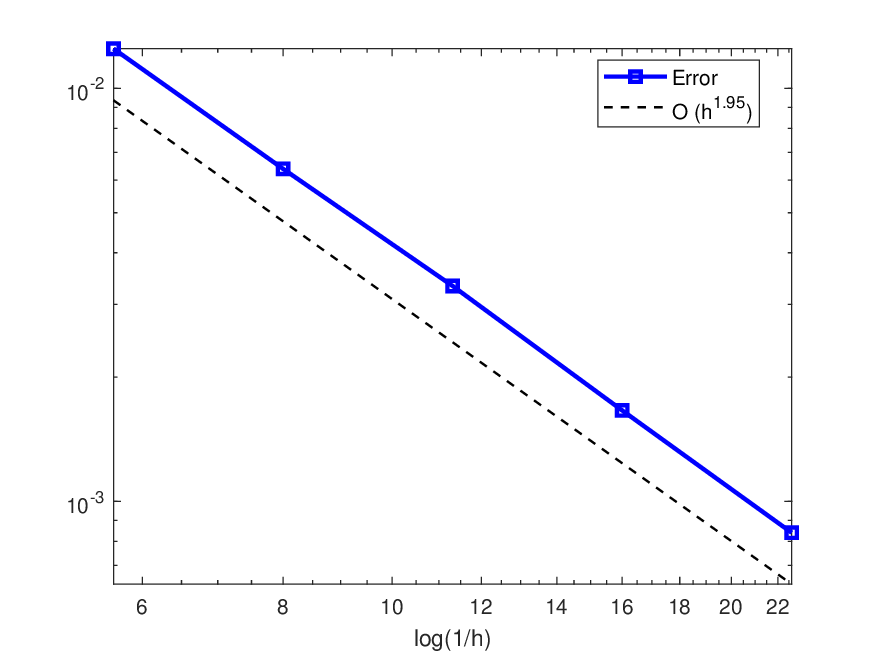}
  \includegraphics[scale=0.35]{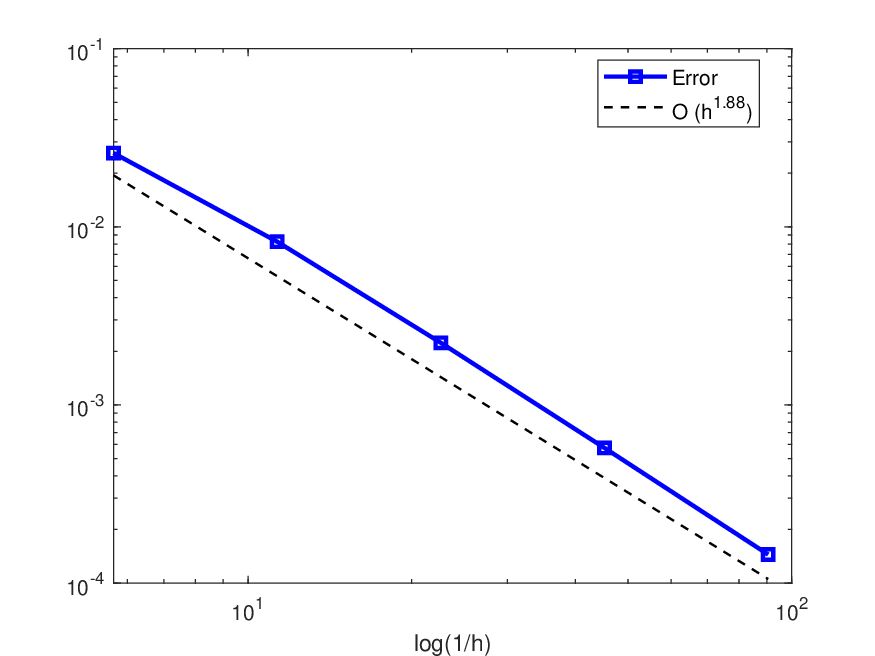}\\
  \caption{The error orders for Example \ref{divfree}. Left: Polygonal meshes; Right: Triangular meshes.}\label{fig:SGEReduced}
\end{figure}

\begin{example}\label{divfree}
This example is given in \cite{CHH2022}. The exact solution of the reduced problem \eqref{reducedProblem} is set to be a divergence-free function in the form
 \[\boldsymbol{u} = \left[ \begin{array}{*{20}{c}}
-x^2 (1-x)^2 y (1-y) (1-2 y) \\
x(1-x)(1-2x)y^2(1-y)^2
\end{array} \right]
.\]
 The right side term $\boldsymbol{f}$ computed from \eqref{reducedProblem} is independent of both $\lambda$ and $\iota$. We use
this $\boldsymbol{f}$ as the right side function of the SGE problem.
\end{example}

In view of the second estimate in Theorem \ref{thm:robust}, we instead compute $E_{u_0} = \|\boldsymbol{u}_0 - \boldsymbol{u}_h \|_{\iota,h}/\|\boldsymbol{f}\|_0$. The convergence orders for $\iota = 10^{-5}$, $\lambda = 10^{8}$ and $\mu = 1$ are shown in Fig.~\ref{fig:SGEReduced} for both the polygonal and triangular meshes. It can be seen that at this time the method is of second-order convergence, although the theoretical prediction in Theorem \ref{thm:robust} only gives $E_{u_0} \le Ch^{1/2}$ when $\iota \le h$. Compared with the finite elements given in \cite{CHH2022}, the virtual element methods exhibit better performance for this example.
Such a phenomenon is also observed for the second example in \cite{Li-Ming-Shi-2018} when our VEM is employed. It should be pointed out that Example \ref{example1} shows that our method only provides a half-order convergence rate in general.
We remark that the half-order convergence rate is not optimal, and the optimal second-order convergence for small $\iota$ can be restored by using the Nitsche's technique.

\section*{Acknowledgments}

JG Huang was partially supported by NSFC (grant 12071289)
and the Fundamental Research Funds for the Central Universities.
Y Yu was partially supported by the National Science Foundation for Young Scientists of China (No. 12301561) and
the key project of the Department of Education of Hunan Province (No. 24A0100).




\bibliographystyle{plain}
\bibliography{RefsD}







\end{document}